\documentclass[12pt]{article}
\usepackage{amsfonts,epsfig,latexsym, multicol, amsmath, vmargin}
\setpapersize{custom}{21cm}{29.7cm}
\setmarginsrb{2cm}{2cm}{2cm}{2cm}{0pt}{0pt}{0pt}{0pt}

\newtheorem{theor}{Theorem}

\newtheorem{lemma}{Lemma}

\newtheorem{conj}{Conjecture}

\begin{document}

\title{Small cones of oriented semi-metrics}

\author{
Michel DEZA\footnote{Michel.Deza@ens.fr} \\
     CNRS/ENS, Paris and Institute of Statistical Mathematics, Tokyo\\
\and Mathieu Dutour \footnote{Mathieu.Dutour@ens.fr}\\
\'Ecole Normale Sup\'erieure, France \\
\and   Elena PANTELEEVA \footnote{pantel@mccme.ru}\\
Moscow Pedagogical University,  Moscow, Russia
}

\maketitle 

\begin{abstract}
We consider polyhedral cones, associated with 
quasi-semi-metrics (oriented distances), in particular, with oriented 
multi-cuts, on $n$ points.
We computed the number of facets and of extreme rays, their adjacencies, and 
 incidences of the cones $QMET_{n}$ and $OMCUT_{n}$ for $ n = 3, 4, 5$ (see Table \ref{tab:smallDatas})
and, partially, for $n=6$. Some results for general $n$ are also given.

\end{abstract}

\section{Introduction and basic notions}
The notions of directed distances, quasi-metrics and oriented multi-cuts are 
 generalizations of the notions of distances, metrics and cuts.
The notions of distances and  metrics are central objects in 
Graph Theory and Combinatorial Optimization. Quasi-metrics are used 
in the Semantics of Computation (see, for example, 
~\cite{Se}) and in Computational Geometry (see, for
example, ~\cite{ACLM}). Oriented distances have been used 
already in \cite{H}, pages 145--146. 

The knowledge of extreme rays of $QMET_n$ 
(cone of all quasi-semi-metrics on $n$ points) for 
small $n$ will allow to build a theory of multi-commodity flows on oriented 
graphs, as well as it was done for non-oriented graphs using dual $MET_n$ 
(cone of all semi-metrics on $n$ points).
 For general theory of quasi-metrics see, for example,~\cite{Fr}, \cite{H},
 ~\cite{W}, and \cite{St}; for general theory of metrics see ~\cite{B},
 ~\cite{DL}.

Set $V_{n}=\{1, \dots, n\}$. Define a {\em quasi-semi-metric} on $V_n$ as a function $d:V_n^{2}\rightarrow R$, with $d_{ii}=0$, satisfying
\begin{itemize}
\item the {\it oriented triangle inequality} $OT_{ij, k}:= d_{ik} + d_{kj} - d_{ij} \geq 0$;
\item the {\it non-negativity inequality} $NN_{ij}:= d_{ij} \geq 0$.
\end{itemize}
If all non-negativity inequalities are strict, then $d$ is called a {\it quasi-metric}. If for all 
$i, j \in V_n^2$, $d_{ij} = d_{ji}$, then $d$ is called a {\it semi-metric}.

The convex cone $QMET_{n}$ is defined by $n(n-1)(n-2)$ oriented triangle inequalities $OT_{ij, k}$ and $n(n-1)$ non-negativity inequalities $NN_{ij}$ on $n$ points; it is a full-dimensional cone in $R^{n(n-1)}$.

Note, that in the semi-metric case, the triangle inequalities imply non-negativity of the distance; it is not the case in our oriented case.

Through all paper we will represent vectors in $R^{n(n-1)}$ as square matrices of order $n$ with zeroes on the main diagonal. It is easy to see that 
\begin{equation*}
MET_n=\{a+a^{T}\mbox{,~with~}a\in QMET_n\}
\end{equation*}
where $a^T$ denotes the transpose matrix. Moreover, any extreme ray $e$ of $MET_n$ has a form $e=g+g^T$, where $g$ is an extreme ray of $QMET_n$ (in fact, $e$ is a sum of extreme rays of $QMET_n$ with non-negative coefficients, and $e=(e+e^T)/2$).
\newline

Consider now the notions of oriented cut and oriented multi-cut quasi-semi-metrics. 
Given an ordered partition  $(S_{1}, \ldots, S_{q})$ ($q \geq 2$) of 
$V_n$, the quasi-semi-metric $\delta^{'}(S_{1}, \ldots, 
S_{q})$ is called an {\it oriented multi-cut}, if 
$\delta^{'}(S_{1}, \ldots, S_{q})_{ij} = 1$ for $i \in S_{\alpha}$, $j 
\in S_{\beta}$, $\alpha < \beta$ and 
$\delta^{'}(S_{1}, \ldots, S_{q})_{ij} = 0$, 
otherwise. This notion was considered, for example, in ~\cite{SL}. Given a subset $S$ of $V_{n}$, the quasi-semi-metric  $\delta^{'}(S)$ is called {\it oriented cut}, if $\delta^{'}(S)_{ij} = 1$ for $i \in S, j \not\in S$
and $\delta^{'} (S)_{ij} = 0$, otherwise. Clearly, an oriented cut is the case $q=2$ of an oriented multi-cut. 
The full-dimensional cone in $R^{n(n-1)}$, generated by all non-zero oriented multi-cuts on $V_{n}$, is denoted by $OMCUT_{n}$.

Given a partition  $(S_{1}, \ldots, S_{q})$ of $V_n$, the semi-metric $\delta(S_{1}, \ldots, S_{q})$ is called a {\it multi-cut}, if 
$\delta(S_{1}, \ldots, S_{q})_{ij} = 1$ for $i \in S_{\alpha}$, $j 
\in S_{\beta}$, $\alpha \not= \beta$ and 
$\delta(S_{1}, \ldots, S_{q})_{ij} = 0$, otherwise. Given a subset $S$ of $V_{n}$, the semi-metric $\delta(S)$ is called a {\it cut}, if $\delta(S)_{ij}=1$ if $|S\cap \{i, j\}|=1$ and $\delta(S)_{ij}=1=0$, otherwise. Clearly, cuts are multi-cuts with $q=2$.
The $n(n-1)/2$-dimensional cone in $R^{n(n-1)}$, generated by all non-zero cuts on $V_{n}$, is denoted by $CUT_{n}$.

Note, that in the semi-metric case, the multi-cuts with $q>2$ are not extreme rays (they are interior points) of $CUT_n$; it relies on the formula:
\begin{equation*}
\delta(S_{1}, \ldots, S_{q})=(\sum_{i=1}^q \delta(S_i))/2
\end{equation*}
On the other hand,
\begin{equation*}
\delta(S_{1}, \ldots, S_{q})=\delta^{'}(S_{1}, \ldots, S_{q})+\delta^{'}(S_{q}, \ldots, S_{1})
\end{equation*}
and so, 
\begin{equation*}
CUT_n=\{a+a^{T}\mbox{,~with~}a\in OMCUT_n\}\;.
\end{equation*}
The number of all oriented cuts on $n$ points is 
$2^{n}$ and the number $p^{'} (n)$ of all oriented multi-cuts on $V_n$
 is the number of all ordered 
partitions of $n$. In fact,  $p'(n) = \frac{1}{2} \sum_{\pi \in Sym(n)} 
x^{1+d(\pi)}$ with $d(\pi) :=|\{i \leq n | a_{i} > a_{i+1} \} |$ 
for the permutation 
$\pi = (a_{1}, \dots , a_{n}) \in Sym(n)$. So, 
$p'(3)=13$, $p'(4)=75$, $p'(5)=541$, $p'(6)=4683$,  $p'(7)=47293$, $p'(8)=545835$. The number of extreme rays of $OMCUT_n$ is $p'(n)-1$, since we excluded $0$; the number of orbits of those rays is $2,5,9,19,35,71$ for $n=3,4,5,6,7,8$.
\newline

Let $C$ be a cone in $\mathbb{R}^{n}$. Given $v \in \mathbb{R}^{n}$, the 
inequality $\sum_{i=1}^nv_ix_i \geq 0$ is said to be {\it valid} for $C$, 
if it holds for all $x \in C$. Then the set 
$\{ x \in C \vert \sum_{i=1}^nv_ix_i = 0 \}$ 
is called the {\it face} of $C$, induced by the valid inequality 
$\sum_{i=1}^nv_ix_i \geq 0$. A face of dimension $\dim(C) - 1$ is called a 
{\it facet} of $C$; a face of dimension $1$ is called an 
{\it extreme ray} of $C$.

Two extreme rays of $C$  are said to be {\it adjacent} on $C$, if
they span a two-dimensional face of $C$.
 Two facets of $C$ are said to be {\it adjacent}, 
if their intersection has dimension $\dim(C) - 2$ (or codimension $2$). 

The incidence number of a facet (or of an extreme ray) is the number of 
extreme rays lying on this facet (or, respectively, of facets containing 
this extreme ray).

The {\it skeleton} graph of the cone $C$ is the graph 
$G_{C}$ whose node-set is the set of  extreme rays of $C$ and with an
 edge between two nodes if they are adjacent on  $C$. 
The {\it ridge} graph of $C$ is the graph $G^{*}_{C}$  whose node-set is
 the set of facets of $C$ and with an edge between two nodes if
 they are adjacent on $C$. So, the ridge graph of a cone  is
 the skeleton of its dual.
\newline

The cones $QMET_n$ and $OMCUT_n$ have the group $Sym(n)$ of all permutations
as a symmetry group. But another symmetry,  called {\it reversal}, exists:
associate to each ray $d$ the ray $d^{r}$ defined by $d^{r}_{ij}=d_{ji}$, 
i.e., in matrix terms, the reversal corresponds to the transposition of 
matrices.
This yields that the group $Z_2\times Sym(n)$ is a symmetry group of the cones $QMET_n$ and $OMCUT_n$. We conjecture that this group is the full symmetry group of those cones; it has been checked by computer for $n=3, 4, 5$. All orbits of faces considered below, are under action of this group.

The {\it representation matrix} of skeleton (or ridge) graph is the 
square matrix where on the place $i,j$ we put the number of members of 
orbit $O_j$ of extreme rays (or facets, respectively), which are adjacent 
to a fixed representative of orbit $O_i$

Comparing with the cones of semi-metrics, the amount of computation and
memory is much bigger in the oriented 
case, because the dimension of the cones $OMCUT_n$ and $QMET_n$ is twice those of $CUT_n$ and $MET_n$, and because oriented multi-cuts with $q>2$ are not interior points of the cone, denoted by $OCUT_n$, generated by oriented cuts. So-called {\it combinatorial explosion} starts from $n=5$ (see Table \ref{tab:smallDatas}) while for corresponding semi-metric cones it starts from $n=7$. All computations were done using the programs {\it cdd} (see~\cite{Fu}) and an adaptation, by the second author, of adjacency decomposition method from \cite{CR}.
\newline 

We consider the cones $OMCUT_n$ and $QMET_n$ for $n=3,4,5$; we give the 
number of facets and of extreme rays for them, their orbits 
 and Tables of their adjacencies and incidences.
We study (having the semi-metric case in mind, see ~\cite{DD4}, ~\cite{DD5}, ~\cite{DDF}, and~\cite{DL}) the skeleton graphs and the ridge graphs of these 
cones, i.e. their diameters, adjacency conditions.

In Table \ref{tab:smallDatas} there is a synthetic summary of the most important information concerning the considered cones; the diameters there are those of the skeleton and ridge graph, respectively.

In Appendix $1$ we represent all $229$ orbits of extreme rays for $QMET_5$ with adjacency, incidence of their representatives and their size. In Appendix $2$ we give the same information for facets of $OMCUT_5$. For representation matrices detailing orbit-wise adjacencies and additional information on those cones see {\it http://www.geomappl.ens.fr/\texttt{\~}dutour/}.

This paper is a follow-up to~\cite{DePant}, which initiates this subject.
\begin{table}
\small
 \begin{tabular}{|c|c|c|c|c|c|} 
 \hline \hline
 cone &dim. &\# of ext. Rays (orbits) &\# of facets (orbits)&diameters \\ 
\hline
 $CUT_3$=$MET_3$ &3&3(1)&3(1)&1; 1 \\ 
 $CUT_4$=$MET_4$ &6&7(2)&12(1)&1; 2 \\ 
 $CUT_5$ &10&15(2)&40(2)&1; 2 \\ 
 $MET_5$ &10&25(3)&30(1)&2; 2 \\ 
 $CUT_6$ &15&31(3)&210(4)&1; 3 \\
 $MET_6$ &15&296(7)&60(1)&2; 2 \\
 $CUT_7$ &21&63(3)&38780(36)&1; 3 \\
 $MET_7$ &21&55226(46)&105(1)&3; 2 \\
 $CUT_8$ &28&127(4)&$\ge 49604520(\ge 2169)$&1; 3 or 4? \\
 $MET_8$ &28&$\ge 119269588, (\ge 3918)$&168(1)&3?; 2 \\
 \hline \hline
 $OMCUT_3$=$QMET_3$ &6&12(2)&12(2)&2; 2 \\
 $OMCUT_4$ &12&74(5)&72(4)&2; 2 \\
 $QMET_4$ &12&164(10)&36(2)&3; 2 \\
 $OMCUT_5$ &20&540(9)&35320(194)&2; 3 \\
 $QMET_5$ &20&43590(229)&80(2)&3; 2\\
 $QHYP_5$ &20&78810(386)&90(3)&4; 2\\
 $OMCUT_6$ &30&4682(19)&$\ge 217847040(\ge 163822)$&2; ?\\
 $QMET_6$ &30&$\ge 182403032(\ge 127779)$&150(2)&?; 2\\
 \hline \hline
\end{tabular}
\caption{Small cones}
\label{tab:smallDatas}
\end{table}

\section{General results about $QMET_n$}
We got by computations the following new facts about small cones $MET_n$:

(i) the full symmetry group of $MET_n$ is $Sym(n)$ for $n=3, 5, ..., 14$, and for $n=4$ it is $Sym(4)\times Sym(3)$,

(ii) the diameter of $G_{MET_7}$ is three,

(iii) \cite{DFPS} obtained a list of $1550825600$ rays of $MET_8$; we computed that this list consists of $3918$ orbits of extreme rays under the symmetry group $Sym(8)$ of this cone,

(iv) the cone $MET_7$ has $46$ orbits of extreme rays and not $41$ as, by a technical mistake, was given in \cite{G} and \cite{DL}.

Let us define {\em oriented $l_{p}$-norm} be: $||x - y||_{p,\,oriented} = \sqrt[p]{\sum_{k=1}^{n} [\max(x_k - y_k,0)]^p}$ and, for $p=\infty$, 
$||x - y||_{\infty,\,oriented} = \max_{k=1}^{n} \max(x_k - y_k,0)$.

\begin{theor}\label{Prop-by-Panteleeva}
(i) Any quasi-semi-metric $d'$ on $n$ points is embeddable in $l^{n}_{\infty,\,oriented}$.

(ii) Any quasi-semi-metric $d'$ on $n$ points is embeddable in $l^m_{1,\,oriented}$ for some $m$ if and only if $d'\in OCUT_n$ (the cone generated by all oriented cuts on $V_n$).
\end{theor}
{\it Proof.} (i) Let $v_1,\dots, v_n$ in $R^{n}$ defined as $v_i = (d'(i, 1), d'(i,2), ..., d'(i,n))$.

Then $||v_i - v_j ||_{\infty,\,oriented} = max(d'(i, k) - d'(j, k),0)$ (it is $\leq d'(i, j)$ from oriented triangle inequality) and 
$d'(i, j) - d'(j, j) = d'(i, j)$ so, $||v_i - v_j ||_{\infty,\,oriented}=d'(i,j)$.

(ii) The proof is the same as in proposition 4.2.2 of \cite{DL} for non-oriented case.

\begin{theor}\label{At-least-n-1-zeros}
Every extreme ray of $QMET_n$ has at least $n-1$ coordinates equal to zero. This lower bound is met for $n=4, 5, 6$.
\end{theor}
{\it Proof.} The rank of the system $(d_{ij}=d_{ik}+d_{kj})_{1\leq i, j, k\leq n}$ is $(n-1)^2$ (see \cite{DePant}). Let $d$ be an extreme ray of $QMET_n$ and let $NN=(NN_{\alpha})_{\alpha\in A}$ be the set of all non-negativity facets, to which $d$ is incident and $OT=(OT_{\beta})_{\beta\in B}$ be the set of all oriented triangle facets, to which it is incident. So, the rank of $NN\cup OT$ is $n(n-1)-1$.

If $rank\,(OT)=(n-1)^2$, then the vector $d$ is incident to all oriented triangle inequalities so, $d_{ij}+d_{ji}=0$, and since $d$ belongs to $QMET_n$, the equalities $d_{ij}=d_{ji}=0$ hold. So, we have $rank\,(OT)\leq (n-1)^2-1$ and
\begin{equation*}
n(n-1)-1=rank\,(NN\cup OT)\leq rank\,(OT)+rank\,(NN)\leq (n-1)^2-1+|A|
\end{equation*}
implying $n-1\leq |A|$.
\\

Above theorem imply that any extreme ray of $QMET_n$ is not the directed path distance of a oriented graph (see, for example \cite{CJTW}, for those notions).

The {\em vertex-splitting} of a vector $(d_{ij})_{1\leq i\not=j\leq n}$ denoted by $(d^{vs}_{ij})_{1\leq i\not= j\leq n+1}$ and defined by
\begin{equation*}
d^{vs}_{n\,n+1}=d^{vs}_{n+1\,n}=0\mbox{,~}d^{vs}_{i\,n+1}=d_{i\,n}
\mbox{~and~}d^{vs}_{n+1\,i}=d_{n\,i}\;.
\end{equation*}
The vertex-splitting of an oriented multicut $d=\alpha(S_1,\dots, S_q)$ is the ray $d^{vs}=\alpha(S_1, \dots, S_l\cup\{n+1\},\dots, S_q)$ if $n\in S_l$. 

\begin{theor}\label{Theorem-vertex-splitting}
The vextex-splitting of an extreme ray of $QMET_n$ is an extreme ray of $QMET_{n+1}$.
\end{theor}
{\it Proof}. If $d$ is an extreme ray of $QMET_n$ then one can check easily the validity of oriented triangle and non-negativity inequalities for $d^{vs}$.\\
The ray $d^{vs}$ is incident to the oriented triangle inequalities $OT_{n,i;n+1}$, $OT_{i,n;n+1}$, $OT_{n+1,i;n}$ and $OT_{i,n+1;n}$ and to the non-negativity inequalities $NN_{n,n+1}$ and $NN_{n+1,n}$.\\
Assume now that $e$ is a ray of $QMET_{n+1}$, which is incident to all facets incident to $d^{vs}$. Then we obtain 
\begin{equation*}
e_{n+1,i}=e_{n,i},\;\; e_{i,n+1}=e_{i,n}\mbox{~and~}e_{n,n+1}=e_{n+1,n}=0
\end{equation*}
Since the restriction of $e$ on $V_n$ yield a function that satisfy to all inequalities, which are satified by $d$, we obtain (since $d$ is an extreme ray) that the restriction $e_{|V_n}$ is a multiple of $d$. So we get $e_{|V_n}=\lambda d$ with $\lambda \geq 0$. Above equalities yield $e=\lambda d^{vs}$.

\begin{theor}
There is no symmetric extreme rays in $QMET_n$
\end{theor}
{\it Proof} If $d$ is an extreme ray then by Theorem \ref{At-least-n-1-zeros} this extreme rays has a zero, say $d_{n\,n+1}=0$. By symmetry one obtain $d_{n+1\,n}=0$, then the validity of triangle inequality $OT_{n,i;n+1}$, $OT_{i,n;n+1}$ give us $d_{i\,n}=d_{i\,n+1}$ and then again by symmetry $d_{n\,i}=d_{n+1\,i}$. So, $d$ is the vertex-splitting of a ray of $QMET_{n-1}$ which will be again a symmetric extreme ray. One can then conclude by induction.

\begin{lemma}
The oriented multicut $\alpha(\{1\},\{2\},\dots, \{n\})$ is an extreme ray of $QMET_n$.
\end{lemma}
{\it Proof}. Let $d=\alpha(\{1\},\{2\},\dots, \{n\})$, we get $d_{ij}=1$ if $j>i$, and $d_{ij}=0$ if $j<i$. The vector $d$ satisfy to all oriented triangle and non-negativity inequalities; so it is a ray of $QMET_n$. Now, let us prove that it is an extreme ray of $QMET_n$. The ray $d$ is incident to:
\begin{enumerate}
\item non-negativity inequalities $NN_{ij}$ with $j<i$
\item oriented triangle inequalities $OT_{i,k;j}$ with $(i,j,k)$ any non-increasing sequence of numbers
\end{enumerate}
Let $e$ be a ray of $QMET_n$, which is incident to the same inequalities of $QMET_n$. We will prove now that $e$ is a non-negative multiple of $d$. Since $e$ is incident to the non-negativity inequalities $NN_{ij}$ with $j<i$, one has $e_{ij}=0$ if $j<i$. It remains to prove that the numbers $e_{ij}$ with $j>i$ are all equal. One has $e$ incident to $OT_{i,j;k}$ with $i<k<j$, so $e_{ik}=e_{ij}$; in the same way $e_{ik}=e_{jk}$ if $j<k$, so we get the result, i.e. $e=\lambda d$ and $\lambda$ is positive since the inequality $e_{12}\geq 0$ is valid.

\begin{theor}
Oriented multi-cuts are extreme rays of $QMET_n$
\end{theor}
{\it Proof.} By repeated applications of above Lemma, Theorem \ref{Theorem-vertex-splitting}, and using group of symmetries, the Theorem follows.
\newline

Two vectors are said to be {\it conflicting} if 
there exist a component on which they have non-zero values
of different sign.

\begin{theor}\label{theorem-easy-cases}
For the ridge graph $G^{*}_{QMET_{n}}$ holds:
    
   (i) A triangle facet is non-adjacent to any facet to which it conflicts;

   (ii) the non-negativity facets $NN_{ij}$ and $NN_{i'j'}$ are non-adjacent if, either $i'=j$, or $j'=i$. 

\end{theor}
{\it Proof.} (i) If $d$ is a quasi-semi-metric vector, which is incident to both, $OT_{ij, k}$ and $NN_{ij}$, then we have $0=d_{ik}+d_{kj}$. But since $d$ is a quasi-semi-metric, we have $d_{ik}\geq 0$ and $d_{kj}\geq 0$; so, $d_{ki}=d_{kj}=0$. The vector $d$ must lie in a space of dimension $n(n-1)-3$. So the facets $OT_{ij, k}$ and $NN_{ij}$ are non-adjacent.\\
If $d$ is a quasi-semi-metric vector incident to $OT_{ij, k}$ and $OT_{ik, l}$ then we will have
\begin{equation*}
d_{ij}=d_{ik}+d_{kj}=d_{il}+d_{lk}+d_{kj}=d_{il}+(d_{lk}+d_{kj})
\end{equation*}
Since $d$ is a quasi-semi-metric we have $d_{ij}\leq d_{il}+d_{lj}$ and $d_{lj}\leq d_{lk}+d_{kj}$. These inequalities must be equalities in order to meet above equality. So, the vector $d$ belongs to a space of dimension $n(n-1)-4$ and the facets are not adjacent.\\
(ii) is obvious, since a vector $d$ incident to $NN_{ij}$ and $NN_{ki}$, is also incident to $NN_{kj}$ and give a lower than expected rank. Similarly, if $d$ satisfy $d_{ij}=0$ and $d_{ji}=0$ then we have the equality $d_{ik}=d_{jk}$ for all $k$ and we will again a too low rank again.

\begin{conj}

   (i) A triangle facet is adjacent to a facet if they are non-conflicting;

   (ii) the non-negativity facets $NN_{ij}$ and $NN_{i'j'}$ are adjacent if, neither $i'=j$, nor $j'=i$;

   (iii) the diameter of $G^*_{QMET_n}$ is $2$.

   (iv) The ridge graph $G^{*}_{QMET_{n}}$ is induced subgraph of $G^{*}_{OMCUT_{n}}$.
\end{conj}
We checked by computer (i)-(iv) for $n\leq 7$. Easy to see that (iii) is 
implied by the criterion of adjacency given by (i) and (ii) together with 
Theorem \ref{theorem-easy-cases}. Apropos, any two facets of 
$QMET_5$ will be not adjacent if, instead of all oriented multi-cuts, 
we restrict ourself only to oriented cuts.

Note, that the set $E_{n} = \{ e + e^{T} | e \mbox{ is an extreme ray of } QMET_{n}\}$ consists, for $n = 3, 4, 5$ of $1, 7, 79$ orbits 
(amongst of $2, 10, 229$), including $0, 3, 10$ orbits of path-metrics of graphs. More 
exactly, a path-metric $d(G)$ belongs to $E_{4}$ for the graphs 
$G=K_{4}, P_{2}, C_{4}, P_{4}$. 

Now, $d(G) \in E_{5}$ for $G = \{K_{2, 3}, K_{5} - K_{3}, K_{5} - P_{2} - P_{3}, 
K_{5}, C_{5}, \overline{P_{2}}, \overline{P_{3}}, \overline{P_{4}}, \overline{P_{5}}, \overline{2P_{2}} \}$, 
where $d(K_{2, 3})$ is an extreme ray of $MET_{5}$;
$d(K_{2, 3})$,  $d(K_{5} - K_{3})$ and $d(K_{5} - P_{2} - P_{3})$ 
do not belong to $CUT_{5}$, and the remaining seven graphs belong 
to $CUT_{5}$. In fact, those seven path metrics $d(G)$ are all of
form $e+e^T$, were $e$ is an extreme ray of $QHYP_5$ (see definition 
of the cone $QHYP_5$ in the end of section five).

\section{General results about $OMCUT_n$}
\begin{conj}

 (i) All oriented multi-cuts are extreme rays of $QMET_n$;

 (ii) the oriented cuts form a {\it dominating clique} (i.e. a complete subgraph such that any node is adjacent to one of its elements) in $G_{OMCUT_n}$; so, the diameter of $G_{OMCUT_n}$ is $2$ or $3$ for any $n\geq 4$.
\end{conj}
We checked it by computer for $n\leq 7$, using only non-negativity facets and oriented triangle facets. The adjacency between oriented cuts and any other oriented multi-cut are preserved even if we consider only oriented triangle and non-ngativity facets, but already for $n=4$ (when $OMCUT_n$ has other facets) the adjacency is dimished for any oriented multi-cut orbit different from oriented cut.
\newline

We have the inclusion $OMCUT_n\subseteq QMET_n$ with equality if and only if $n=3$. 

Let us consider the following inequalities

\begin{itemize}
\item the {\it zero-extension of an inequality} $\sum_{1\leq i\not= j\leq n-1} f_{ij}d_{ij}\geq 0$, is an inequality
\begin{equation*}
\sum_{1\leq i\not= j\leq n}f'_{ij}d_{ij}\geq 0\mbox{~with~} f'_{ni}=f'_{in}=0\mbox{~and~} f'_{ij}=f_{ij}\mbox{,~otherwise};
\end{equation*}
\item the inequality $A_{n}(c_{1},\dots, c_{n-2}; a,b):= \sum^{n-2}_{i=1} ( d_{ac_{i}} + d_{c_{i}b}) + d_{ba} \geq
S_{ab} + S_{c_{1} \dots c_{n-2}}$, where $S_{c_1, \dots c_{n-2}}$ 
denotes the sum of distances along oriented cycle $c_1, \dots ,c_{n-2}$;
\item the inequality $B_{n}(c_{1},\dots, c_{n-2}; a,b):= \sum^{n-2}_{i=1} (d_{c_{i}a} + d_{ac_{i}} 
+ d_{c_{i}b}) \geq d_{ab} + S_{c_{1} \dots c_{n-2}}$;
\item the {\it hypermetric} inequality $H(b):= \sum_{1 \leq i \not=   j \leq n} b_{i}b_{j}d_{ij} \leq 0$, where $b=(b_{1}, \dots , b_{n}) 
\in \mathbb{Z}^{n}$, $\sum^{n}_{i=1} b_{i} = 1$.
\end{itemize}

Call a face of $OMCUT_n$ {\em symmetric} if, in matrix terms, it is preserved by the transposition.

\begin{theor}
The following inequalities are valid on (i.e. faces of) $OMCUT_n$:
     
(i) zero-extensions of valid faces of $OMCUT_{n-1}$;

(ii) symmetric faces coming from $CUT_n$ (so, including any inequality $H(b)$)

(iii) $A_n(c_1, \dots, c_n; a, b)$ and $B_n(c_1, \dots, c_n; a, b)$.
\end{theor}
{\it Proof} (ii) Since we have $CUT_n=\{d+d^r\mbox{~s.t.~}d\in OMCUT_n\}$, the validity of symmetric faces coming from $CUT_n$ is obvious.

(iii) Consider the inequality $A_{n}$ on an o-multi-cut 
$\delta^{'}(S_{1}, \dots, S_{t})$. In this case 
$S_{c_{1} \dots c_{n-2}} \leq n - 3$ and $S_{c_{1} \dots c_{n-2}} = n - 3$ 
if and only if 
$c_{1} \in S_{\alpha_{1}}, \dots, c_{n-2} \in S_{\alpha_{n-2}}$,
 where $1 \leq \alpha_{1} < \dots < \alpha_{n-2} \leq n-2$. 

Let $b \in S_{\alpha}, a \in S_{\beta}$, where $\alpha < \beta$. 
Then $S_{ab} = d_{ba} = 1$ and $S_{ab} + S_{c_{1} \dots c_{n-2}} \leq n-2$. 
But 
$\sum^{n-2}_{i=1} (d_{ac_{i}} + d_{c_{i}b} + d_{ba}) \geq \sum^{n-2}_{i=1} d_{ba} = n-2$ 
and $A_{n}$ holds.

Let now $a \in S_{\alpha}, b \in S_{\beta}$, where $\alpha < \beta$. Then 
$S_{ab} = d_{ab} = 1$ and $S_{ab} + S_{c_{1} \dots c_{n-2}} \leq n-2$. 
But in 
this case $d_{ac_{i}} + d_{c_{i}b}  \geq 1$, hence 
$\sum^{n-2}_{i=1} (d_{ac_{i}} + d_{c_{i}b} + d_{ba}) = \sum^{n-2}_{i=1} 
(d_{ac_{i}} + d_{c_{i}b}) \geq n-2$,   and $A_{n}$ holds.

If $a, b \in S_{\alpha}$, then $S_{ab} = 0$ and 
$S_{ab} + S_{c_{1} \dots c_{n-2}} \leq n-3$. But in this case   
$d_{ac_{i}} + d_{c_{i}b} = 1$  for $c_{i} \notin S_{\alpha}$ and 
$d_{ac_{i}} + d_{c_{i}b} = 0$  for $c_{i} \in S_{\alpha}$, hence, 
$\sum^{n-2}_{i=1} (d_{ac_{i}} + d_{c_{i}b} + d_{ba}) = \sum^{n-2}_{i=1} 
(d_{ac_{i}} + d_{c_{i}b}) \geq S_{c_{1} \dots c_{n-2}}$ and $A_{n}$ holds.

Consider the inequality $B_{n}$ on an o-multi-cut 
$\delta^{'}(S_{1}, \dots , S_{t})$. In this case 
$S_{c_{1} \dots c_{n-2}} \leq n - 3$ and 
$S_{c_{1} \dots c_{n-2}} \leq n - 3$ if and only if 
$c_{1} \in S_{\alpha_{1}}, \dots, c_{n-2} \in S_{\alpha_{n-2}}$, 
where $1 \leq \alpha_{1} < \dots < \alpha_{n-2} \leq t$. 

Let $b \in S_{\alpha}, a \in S_{\beta}$, where $\alpha < \beta$. 
Then $d_{ab} = 0$ and $d_{ab} + S_{c_{1} \dots c_{n-2}} \leq n - 3$. 
But in this case  
$d_{ac_{i}} + d_{c_{i}b}  \geq 1$, hence 
$\sum^{n-2}_{i=1} (d_{ac_{i}} + d_{c_{i}a} + d_{c_{i}b}) \geq  
\sum^{n-2}_{i=1} 
(d_{ac_{i}} + d_{c_{i}b}) \geq n-2$,   and $B_{n}$ holds.

Let $a \in S_{\alpha}, b \in S_{\beta}$, where $\alpha < \beta$. 
Then $d_{ab} = 1$ and $d_{ab} + S_{c_{1} \dots c_{n-2}} \leq n - 2$. 
In this case  
$d_{ac_{i}} + d_{c_{i}b}  \geq 1$, hence 
$\sum^{n-2}_{i=1} (d_{ac_{i}} + d_{c_{i}a} + d_{c_{i}b}) \geq  
\sum^{n-2}_{i=1} 
(d_{ac_{i}} + d_{c_{i}b}) \geq n-2$,   and $B_{n}$ holds.

If $a, b \in S_{\alpha}$, then $d_{ab} = 0$ and 
$d_{ab} + S_{c_{1} \dots c_{n-2}} \leq n-3$. But in this case   
$d_{ac_{i}} + d_{c_{i}a} +  d_{c_{i}b} \geq  1$  for 
$c_{i} \notin S_{\alpha}$ and $d_{ac_{i}} + d_{c_{i}a} 
+ d_{c_{i}b} = 0$  for $c_{i} \in S_{\alpha}$, hence, 
$\sum^{n-2}_{i=1} (d_{ac_{i}} + d_{c_{i}b} + d_{c_{i}b})
 \geq S_{c_{1} \dots c_{n-2}}$ and $B_{n}$ holds.

\begin{conj}
The following inequalities correspond to facets of $OMCUT_n$:

(i) zero-extensions of any facet of $OMCUT_{n-1}$ (so, including any oriented triangle and non-negativity inequalities);

(ii) any hypermetric inequality $H(b)$, except of non-oriented triangle inequalities;

(iii) $A_{n}(c_1, \dots, c_{n-2}; a,b)$ and $B_{n}(c_1, \dots, c_{n-2}; a,b)$.
\end{conj}
We checked this conjecture by computer for $n\leq 7$.

Any oriented triangle inequality is a zero extension of $A_3$ while any non-negativity inequality is a zero-extension of a degenerated $B_{3}$ with $b=c_1$.
\newline

The first symmetric facet, $H(1, 1, 1, -1, -1)$, is the only symmetric facet of $OMCUT_5$.
\begin{theor}

(i) Any symmetric facet of the cone $OMCUT_n$ correspond to a facet of the cone $CUT_n$.

(ii) All orbits of symmetric facets of $OMCUT_{n}$, $n\leq 6$
are represented by $H(b)$ with $b= (1, 1, 1, -1, -1)$, $(1, 1, 1, -1, -1, 0)$, $(2, 1, 1, -1, -1, -1)$, and $(1, 1, 1, 1, -1, -2)$.

(iii) All orbits of symmetric facets of $OMCUT_{7}$, 
are all $9$ orbits of hypermetric facets of $CUT_7$ different from the orbit of triangle facets, and $18$ out of $26$ non-hypermetric ones (namely, all but $O_6$, $O_{13}$, $O_{22}$, $O_{18}$, $O_{20}$, $O_{24}$, $O_{25}$, and $O_{26}$ in terms of \cite{DD})

\end{theor}
{\it Proof}: (ii) and (iii) were obtained by computer using (i). In order to prove (i), let us fix a symmetric facet $F$ of $OMCUT_n$. We set $U_F=\{d\in OMCUT_n\mbox{~s.t.~}F(d)=0\}$.
If $F$ is a symmetric facet, then $U_F$ is a set invariant by the reversal (transposition).\\
Denote $SU(X):=\{d+d^{r}\mbox{~s.t.~}d\in X\}$ for any $X\subset OMCUT_n$. Then $SU(U_F)$ is a set of semi-metrics, which are incident to $F$. Moreover, since 
$CUT_n=SU(OMCUT_n)$, we have $SU(U_F)\subset OMCUT_n$
By hypothesis $F$ is a facet, so $U_F$ has dimension $n(n-1)-1$. The mapping $d\mapsto d+d^{r}$ decrease dimension by at most $n(n-1)/2$; so, we get that $SU(U_F)$ has dimension $n(n-1)/2-1$, i.e. $F$ is a facet of $CUT_n$.

Given two oriented partitions, $A$ and $B$, we will write $A<B$, if $A$ is a proper refinement of $B$. We will write $Q\dot{<}B$ if, moreover, each part of $A$ is a {\it proper} subset of a part of $B$.

\begin{conj}

(i) An oriented multi-cut $\delta'(A)$ is not adjacent to all oriented multi-cuts $\delta'(B)$ such that $B<A^{T}$.

(ii) The orbit of extreme rays represented by oriented cut $\delta'(\{1\}, \{2, \dots, n\})$ is unique orbit, such that extreme rays in this orbit is not adjacent {\it only} to oriented multi-cuts described in (i) above; the total adjacency is $p'(n)-p'(n-1)-1$ and it is the maximal total adjacency.

(iii) An extreme ray of the orbit, represented by oriented cut $\delta'(\{1, 2\}, \{3, \dots, n\})$, is not adjacent only to oriented multi-cuts $\delta'(B)$ such that either $B< (\{1, 2\}, \{3, \dots, n\})^{T}$, or $B$ is any cyclic shift of $C$ with $C\dot{<} (\{1, 2\}, \{3, \dots, n\})$.

(iv) The diameter of $G_{OMCUT_n}$ is $2$.


\end{conj}

\section{The cases of $3,4$ points}

 We start with complete description  of the cone 
$QMET_{3} = OMCUT _{3}$.

 There are $12$ extreme rays in $OMCUT_{3}$: $12$ non-zero 
oriented multi-cuts on $V_3$, including $6$ oriented cuts.
Under the group action we have, in fact, only $2$ orbits to consider:
 the orbit $O_{1}$ of oriented cuts and the orbit $O_{2}$ of other 
 oriented multi-cuts.
The list of representatives of the orbits 
 is given in the Table
\ref{tab:tabl3O}.

\begin{table}
\begin{center}
\scriptsize
\begin{tabular}{|c|cccccc|}
 \hline
 &12&13&21&23&31&32\\
 \hline 
 $\delta^{'}(\{1\})$:&1&1&0&0&0&0\\
 $\delta^{'}(\{1\}, \{2\}, \{3\})$:&1&1&0&1&0&0\\
 \hline
 $OT_{12, 3}$&-1&1&0&0&0&1\\
 $NN_{12}$& 1&0&0&0&0&0\\
 \hline
\end{tabular}
\caption{Representative of the orbits of extreme rays and facets in $OMCUT_{3}$}
\label{tab:tabl3O}
\end{center}
\end{table}

Note that all oriented cuts above can be obtained from $\delta^{'}(\{1\})$ by 
a permutation ($\delta^{'}(\{2\})$ and $\delta^{'}(\{3\})$) or by a reversal
 and a permutation ($\delta^{'}(\{1, 2\})$, $\delta^{'}(\{1, 3\})$, and 
$\delta^{'}(\{2, 3\})$); all oriented multi-cuts with $q=3$ above can be 
 obtained from $\delta^{'}(\{1\}, \{2\}, \{3\})$ by some permutation.

 The only facet-defining inequalities of $OMCUT_{3}$ are 
 $6$ oriented triangle inequalities $OT_{ij, k}$ and $6$ 
non-negativity inequalities $NN_{ij}$ which form two orbits; see 
their representatives in Table \ref{tab:tabl3O}.

Adjacencies of facets and extreme rays of $OMCUT_{3}$ are shown in Table
 \ref{tab:tabl4}. For each orbit a representative and 
a number of adjacent ones from other orbits are given, as well as the total
 number of adjacent ones, the number of incident extreme rays (respectively, 
facets) and the size of orbits. 
\newline

\begin{table}
\scriptsize
\begin{tabular}{|c|c|cc|c|c|c||c|c|cc|c|c|c|}
 \hline 
Orbit & Representative &  $O_{1}$ & $O_{2}$ & Adj. & Inc. & Orbit size& Orbit & Representative & $OT$ &  $NN$ &  Adj. & Inc. & Orbit size\\ 
 \hline
 $O_{1}$ &  $\delta^{'}(\{3\})$ & $5$ & $4$ & $9$ & $8$ & $6$& $OT$ & $OT_{12, 3}$ & $3$ & $5$ & $8$ & $7$ & $6$\\ 
 $O_{2}$ & $\delta^{'}(\{3\},  \{2\}, \{1\})$ & $4$ & $2$ & $6$ & $6$ & $6$& $NN$ & $NN_{12}$ & $5$ & $2$ & $7$ & $6$ & $7$\\
 \hline 
\end{tabular}
\caption{Representation matrix of $G_{OMCUT_{3}}$ and $G^*_{OMCUT_{3}}$}
\label{tab:tabl4}
\end{table}
 
For the next case $n=4$ we have
\begin{enumerate}
\item $OMCUT_4$ has $74$ extreme rays (all non-zero 
 oriented multi-cuts on $V_{4}$). 
Under the group action we have, in fact, only $5$ orbits to consider. They are orbits  with the representatives 
 $\delta^{'}(\{4\})$ (orbit $O_{1}$), $\delta^{'}(\{4,3\})$ (orbit 
 $O_{2}$) -- the only two orbits of oriented cuts in $OMCUT_4$, 
$\delta^{'}(\{4\}, \{3\}, \{2, 1\})$ (orbit $O_{3}$), 
 $\delta^{'}(\{4\}, \{3, 2\}, \{1\})$ (orbit $O_{4}$) and 
 $\delta^{'}(\{4\}, \{3\}, \{2\}, \{1\})$ (orbit $O_{5}$).
We show these representatives in the Table \ref{tab:tabl4F}.
\item $OMCUT_4$ has $72$ facets from $4$ orbits, which are induced by $24$ oriented triangle inequalities $OT_{ij, k}$, $12$ non-negativity inequalities $NN_{ij}$, $12$ inequalities $A_4(i_1, i_2; j_1, j_2)$, and $24$ inequalities $B_4(i_1, i_2; j_1, j_2)$. See Table \ref{tab:tabl4F} for representatives of the orbits and Table \ref{tab:tabl11} for the representation matrix.
\item The cone $QMET_4$ has $36$ facets, distributed in two orbits: $24$ 
oriented triangle facets (orbit $OT$) and $12$ non-negativity 
facets  (orbit $NN$).
\item There are $164$ extreme rays in $QMET_4$. Under the group action we have $10$ orbits to consider: 
 orbits $O_{1}$ -- $O_{5}$ with the same representatives as 
 in $OMCUT_4$  and $5$ other orbits. The list of representatives of the orbits  is given in the Table \ref{tab:tabl4QO}.
\end{enumerate}

\begin{table}
\begin{center}
\scriptsize
\begin{tabular}{|c|c|cccccccccccc|}
 \hline
Or. &Representative &12&13&14&21&23&24&31&32&34&41&42&43\\
 \hline
 $O_1$ & $\delta'(\{4\})$&0&0&0&0&0&0&0&0&0&1&1&1\\
 $O_2$ & $\delta'(\{3, 4\}, \{1, 2\})$&0&0&0&0&0&0&1&1&0&1&1&0\\
 $O_3$ & $\delta'(\{4\}, \{3\}, \{1, 2\})$&0&0&0&0&0&0&1&1&0&1&1&1\\
 $O_4$ & $\delta'(\{4\}, \{2, 3\}, \{1\})$&0&0&0&1&0&0&1&0&0&1&1&1\\
 $O_5$ & $\delta'(\{4\}, \{3\}, \{2\}, \{1\})$&0&0&0&1&0&0&1&1&0&1&1&1\\
 \hline
 $OT$ & $OT_{12, 3}$&-1&1&0&0&0&0&0&1&0&0&0&0\\
 $NN$ & $NN_{12}$&1&0&0&0&0&0&0&0&0&0&0&0\\
 $A_4$ &$A_4(1,2;3,4)$&-1&0&1&-1&0&1&1&1&-1&0&0&1\\
 $B_4$ & $B_4(1,2;3,4)$&-1&1&1&-1&1&1&1&1&-1&0&0&0\\
 \hline
\end{tabular}
\caption{The representatives of orbits of extreme rays and facets in $OMCUT_4$}
\label{tab:tabl4F}
\end{center}
\end{table}

\begin{table}
\begin{center}
\scriptsize
\begin{tabular}{|c|c|cccc|c|c|c|} 
 \hline 
 Orbit & Representative & $OT$ & $NN$ & $A_4$ & $B_4$ & Adj. & Inc. & Size\\ 
 \hline 
 $OT$ & $OT_{12, 3}$ & $17$ &  $11$ & $5$ & $8$ & $41$ & $43$ & $24$\\
 $NN$ & $NN_{12}$ & $22$ & $6$ & $12$ & $8$ & $48$ & $43$ & $12$\\
 $A_4$ & $A_4(1, 2; 3, 4)$ & $10$ & $12$ & $0$ & $2$ & $24$ & $28$ & $12$\\
 $B_4$ & $B_4(1, 2; 3, 4)$ & $8$ & $4$ & $1$ & $3$ & $16$ & $17$ & $24$\\
 \hline
\end{tabular}
\caption{Representation matrix of $G^*_{OMCUT_4}$}
\label{tab:tabl11}
\end{center}
\end{table}
 
\begin{table}
\begin{center}
\scriptsize
\begin{tabular}{|c|c|ccccc|c|c|c|} 
 \hline 
 Orbit & Representative &  $O_{1}$ & $O_{2}$ & $O_{3}$ & $O_{4}$ & $O_{5}$ & Adj. & Inc. & $\vert  O_{i} \vert $\\ 
 \hline 
 $O_{1}$ & $\delta^{'}(\{4\})$ & $7$ & $6$ & $21$ &  $9$ & $18$ & $61$ & $42$ & $8$\\
 $O_{2}$ & $\delta^{'}(\{4, 3\})$ & $8$ & $5$ & $20$ & $12$ & $8$ & $53$ & $48$ & $6$\\
 $O_{3}$ & $\delta^{'}(\{4\}, \{3\}, \{2, 1\})$ & $7$ & $5$ & $15$ & $7$ & $10$ & $44$ & $29$ & $24$\\
 $O_{4}$ & $\delta^{'}(\{4\}, \{3, 2\}, \{4\})$ & $6$ & $6$ & $14$ & $6$ & $8$ & $40$ & $33$ & $12$\\
 $O_{5}$ & $\delta^{'}(\{4\}, \{3\}, \{2\}, \{1\})$ & $6$ & $2$ & $10$ & $4$ & $12$ & $34$ & $24$ & $24$\\
 \hline 
\end{tabular}
\caption{Representation matrix of $G_{OMCUT_4}$}
\label{tab:tabl12}
\end{center}
\end{table}

Denote by $e(i, j)$ the $\{0, 1\}$-vector in $\mathbb{R}^{12}$
with $1$ only on the place $ij$. In Table \ref{tab:tabl4QO} there is the list of $10$ orbits of extreme rays and $2$ orbits of facets of $QMET_4$ .
\begin{table}
\begin{center}
\scriptsize
\begin{tabular}{|c|c|cccccccccccc|}
 \hline
 Or.&Representative &12&13&14&21&23&24&31&32&34&41&42&43\\
 \hline 
 $O_{1}$ & $\delta^{'}(\{1\})$&0&0&0&0&0&0&0&0&0&1&1&1\\
 $O_{2}$ & $\delta^{'}(\{1, 2\})$&0&0&0&0&0&0&1&1&0&1&1&0\\
 $O_{3}$ & $\delta^{'}(\{1\}, \{2\}, \{3,  4\})$&0&0&0&0&0&0&1&1&0&1&1&1\\
 $O_{4}$ & $\delta^{'}(\{1\}, \{2, 3\}, \{4\})$&0&0&0&1&0&0&1&0&0&1&1&1\\
 $O_{5}$ & $ \delta^{'}(\{1\}, \{2\}, \{3\}, \{4\})$&0&0&0&1&0&0&1&1&0&1&1&1\\
 $O_{6}$ & $\delta'(\{1\}, \{2\}, \{3\}, \{4\})+e(1, 4)$&0&0&0&1&0&0&1&1&0&2&1&1\\
 $O_{7}$ & $\delta'(\{1\}, \{2\}, \{3\}, \{4\})+e(4, 3)$&0&0&0&1&0&0&1&1&1&1&1&1\\
 $O_{8}$ & $\delta'(\{1\}, \{2\}, \{3\}, \{4\})+e(3,2)$&0&0&0&1&0&1&1&1&1&1&1&0\\
 $O_{9}$ & $\delta'(\{1\}, \{2\}, \{3\}, \{4\})+e(2, 1)+e(4, 3)$&0&0&1&1&1&1&1&1&1&1&0&0\\
 $O_{10}$ & $\delta'(\{1\}, \{2\}, \{3\}, \{4\})+e(1,4)+e(2,1)+e(4,3)$&0&0&1&1&1&1&1&1&2&1&0&0\\
 \hline 
 $OT$ & $OT_{12, 3}$&-1&1&0&0&0&0&0&1&0&0&0&0\\
 $NN$ & $NN_{12}$&1&0&0&0&0&0&0&0&0&0&0&0\\
 \hline 
\end{tabular}
\caption{Representatives of orbits of extreme rays and facets in $QMET_4$}
\label{tab:tabl4QO}
\end{center}
\end{table}






 The adjacencies and incidences of the facets and extreme rays of 
 $QMET_4$ are given in Tables \ref{tab:tabl15} -- \ref{tab:tabl16}.

\begin{table} 
\begin{center}
\scriptsize
\begin{tabular}{|c|c|cc|c|c|c|} 
 \hline 
 Orbit & Representative & $OT$ &  $NN$ & Adj. & Inc. & Size\\
 \hline 
 $OT$ & $OT_{12, 3}$ & $17$ & $11$ & $28$ & $78$ & $24$\\
 $NN$ & $NN_{12}$ & $22$ & $6$ & $28$ & $80$ & $12$\\
 \hline 
\end{tabular}
\caption{The adjacencies of facets in $QMET_4$}
\label{tab:tabl15}
\end{center}
\end{table}

\begin{table}
\scriptsize
\begin{tabular}{|c|c|cccccccccc|c|c|c|} 
 \hline 
 Or. &  Representative & $O_{1}$ & $O_{2}$ & $O_{3}$ & $O_{4}$ & $O_{5}$ & 
 $O_{6}$ & $O_{7}$ & $O_{8}$ & $O_{9}$ & $O_{10}$ & Adj. & Inc. & $\vert O_{i} \vert $\\
 \hline
 $O_{1}$&$\delta^{'}(\{1\})$&$7$&$6$&$21$&$9$&$18$&$6$&$9$&$6$&$3$&$6$&$91$&$27$&$8$\\
 $O_{2}$&$\delta^{'}(\{1, 2\})$&$8$&$5$&$20$&$12$&$8$&$12$&$16$& $4$&$4$&$8$&$97$&$24$&$6$\\
 $O_{3}$&$\delta^{'}(\{1\}, \{2\}, \{3,  4\})$&$7$&$5$&$7$&$5$&$10$&$4$&$4$&$2$&$0$&$2$&$46$&$21$&$24$\\
 $O_{4}$&$\delta^{'}(\{1\}, \{2, 3\}, \{4\})$&$6$&$6$&$10$&$2$&$8$&$4$&$4$&$0$&$2$&$4$&$46$&$21$&$12$\\
 $O_{5}$& $ \delta^{'}(\{1\}, \{2\}, \{3\}, \{4\})$&$6$&$2$&$10$&$4$&$3$& $4$&$2$&$1$&$0$&$1$&$33$&$18$&$24$\\
 $O_{6}$&$\delta'(\{1\}, \{2\}, \{3\}, \{4\})+e(1, 4)$&$2$&$3$&$4$&$2$&$4$&$0$&$2$&$1$&$0$&$0$&$18$&$16$&$24$\\
 $O_{7}$&$\delta'(\{1\}, \{2\}, \{3\}, \{4\})+e(4, 3)$&$3$&$4$&$4$&$2$&$2$&$2$&$0$&$1$&$1$&$2$&$21$&$15$&$24$\\
 $O_{8}$&$\delta'(\{1\}, \{2\}, \{3\}, \{4\})+e(3,2)$&$4$&$2$&$4$&$0$&$2$&$2$&$2$&$0$&$0$&$0$&$16$&$15$&$12$\\
 $O_{9}$&$\delta'(\{1\}, \{2\}, \{3\}, \{4\})+e(2, 1)+e(4, 3)$&$4$&$4$&$0$&$4$&$0$&$0$& $4$&$0$&$0$&$4$&$20$&$12$&$6$\\
 $O_{10}$&$\delta'(\{1\}, \{2\}, \{3\}, \{4\})+e(1,4)+e(2,1)+e(4,3)$&$2$&$2$&$2$&$2$&$1$&$0$&$2$&$0$&$1$& $0$&$12$&$12$&$24$\\
 \hline
\end{tabular}
\caption{The adjacencies of extreme rays in $QMET_4$}
\label{tab:tabl16}
\end{table}

 Note, that in $QMET_4$  the adjacencies of facets $OT_{ij, k}$ and 
 $NN_{ij}$ are the same as in $OMCUT_4$ (see Tables \ref{tab:tabl11}
and \ref{tab:tabl15}); hence, $G^{*}_{QMET_4}$  is an induced subgraph of 
$G^{*}_{OMCUT_4}$. But the adjacencies of oriented multi-cuts from orbits $O_{3}$,
  $O_{4}$ and $O_{5}$ are decreased in the cone  $QMET_4$ (see Tables
 \ref{tab:tabl12} and \ref{tab:tabl16}); hence, $G_{OMCUT_4}$ is not an 
induced subgraph of $G_{QMET_4}$.

\section {The case of $5$ points} 

We present here the complete description of $QMET_5$ and of $OMCUT_5$.

 The quasi-semi-metric cone $QMET_5$ has $80$ facets, distributed in two
 orbits: $60$ triangle facets (orbit $F_{1}$) and $20$ non-negativity facets 
 (orbit $F_{2}$).
The list of representatives of the orbits 
 is given in the Table
\ref{tab:tabl5QF}.

\begin{table}
\scriptsize
\begin{tabular}{|c|p{8pt}p{4pt}p{4pt}p{4pt}p{4pt}p{4pt}p{4pt}p{4pt}p{4pt}p{4pt}p{4pt}p{4pt}p{4pt}p{4pt}p{4pt}p{4pt}p{4pt}p{4pt}p{4pt}p{4pt}||c|cc|c|c|c|}
 \hline
 &12&13&14&15&21&23&24&25&31&32&34&35&41&42&43&45&51&52&53&54& Orbit & $OT$ & $NN$ & Adj. & Inc. & $\vert F_{i} \vert$\\
 \hline
 $OT_{1,2; 3}$&-1&1&0&0&0&0&0&0&0&1&0&0&0&0&0&0&0&0&0&0& $OT$ & 19 & 49 & 68 & 13590 & 60\\
 $NN_{12}$&1&0&0&0&0&0&0&0&0&0&0&0&0&0&0&0&0&0&0&0& $NN$ & 12 & 57 & 69 & 14359 & 20\\
 \hline
\end{tabular}
\caption{The list of the orbits of facets in $QMET_5$ and their representation matrix}
\label{tab:tabl5QF}
\end{table}


There are $43590$ extreme rays in $QMET_5$. Under the group action we have  
$229$ orbits to consider (see Appendix $1$ for the full list of representatives); amongst them, nine orbits of oriented multi-cuts on $V_{5}$ are given in the Table \ref{tab:tabl5O}.

Note, that, for example, the representative of orbit $O_{18}$ is similar to the representative
 $\delta^{'}(\{1\}, \{2\}, \{3\}, \{4\}) + e(1, 4)$ 
of the orbit  $O_{6}$ in the cone $QMET_4$: it is 
$\delta^{'}(\{1\}, \{2\}, \{3\}, \{4\}, \{5\}) + e(1, 5)$ , 
and the representative of orbit $O_{11}$ is similar to the representative
 $\delta^{'}(\{1\}, \{2\}, \{3\}, \{4\}) + e(4, 3)$ 
of the orbit  $O_{7}$ in the cone $QMET_4$: it is 
$\delta^{'}(\{1\}, \{2\}, \{3\}, \{4\}, \{5\}) + e(5, 4)$.

Some other representatives from non-oriented multi-cut's orbits of the cone $QMET_5$ have similar form. So, if $v_{i}$ is the representative of orbit $O_{i}$, then

$v_{11} =
\delta^{'}(\{1\}, \{2\}, \{3\}, \{4\}, \{5\}) + e(5, 4)$.

$v_{18} =
\delta^{'}(\{1\}, \{2\}, \{3\}, \{4\}, \{5\}) + e(1, 5)$ ,

 $v_{25} =
\delta^{'}(\{1\}, \{2\}, \{3\}, \{4\}, \{5\}) + e(5, 3) + e (5, 4)$, 

 $v_{26} =
\delta^{'}(\{1\}, \{2\}, \{3\}, \{4\}, \{5\}) +e(1, 5) + e(5, 4)$, 

 $v_{32} =
\delta^{'}(\{1\}, \{2\}, \{3\}, \{4\}, \{5\}) + e(1, 5) + e(2, 5)$, 


 $v_{50} =
\delta^{'}(\{1\}, \{2\}, \{3\}, \{4\}, \{5\}) +  e(1, 4)  +  e(1, 5) +  e(5, 4)$, 

 $v_{52} =
\delta^{'}(\{1\}, \{2\}, \{3\}, \{4\}, \{5\}) + e(5, 3)  + e(5, 4)  + e(4, 3)$, 

 $v_{53} =
\delta^{'}(\{1\}, \{2\}, \{3\}, \{4\}, \{5\}) + e(1, 4)  + e(5, 3)  + e(5, 4)$, 
 



 $v_{84} =
\delta^{'}(\{1\}, \{2\}, \{3\}, \{4\}, \{5\}) + e(1, 4)  + e(1, 5)  + e(2, 5)$,



%
%
%
%
%
%
%

The adjacencies of the facets of $QMET_5$ are  given in Table \ref{tab:tabl5QF}.

Oriented cuts (orbits $O_{1}$ and $O_{2}$ together) form a  clique, but, distinctly from cases $n=3, 4$, not a dominating clique. For
 example, a representative from orbit $O_{108}$ (and the same from  $O_{155}$, $O_{157}$,  $O_{172}$, $O_{185}$, $O_{186}$, $O_{207}$, $O_{216}$ -- $O_{227}$, $O_{229}$)
is not adjacent with any oriented cut.

The cone $OMCUT_5$ has $540$ extreme rays (all non-zero oriented multi-cuts 
on $V_{5}$), which form $9$ orbits.
See Table \ref{tab:tabl5O} for representatives of the orbits and Table \ref{tab:tabl5M} for the representation matrix.

\begin{table}
\scriptsize
\begin{tabular}{|c|cccccccccccccccccccc|}
 \hline
 &12&13&14&15&21&23&24&25&31&32&34&35&41&42&43&45&51&52&53&54\\
 \hline 
 $\delta^{'}(\{1\})$&1&1&1&1&0&0&0&0&0&0&0&0&0&0&0&0&0&0&0&0\\
 $\delta^{'}(\{1,5\})$&1&1&1&0&0&0&0&0&0&0&0&0&0&0&0&0&0&1&1&1\\
 $\delta^{'}(\{1\},\{2\},\{3,4,5\})$&1&1&1&1&0&1&1&1&0&0&0&0&0&0&0&0&0&0&0&0\\
 $\delta^{'}(\{1,2\},\{3\},\{4,5\})$&0&1&1&1&0&1&1&1&0&0&1&1&0&0&0&0&0&0&0&0\\
 $\delta^{'}(\{1\},\{2,3,4\},\{5\})$&1&1&1&1&0&0&0&1&0&0&0&1&0&0&0&1&0&0&0&0\\
 $\delta^{'}(\{1\},\{2,3\},\{4,5\})$&1&1&1&1&0&0&1&1&0&0&1&1&0&0&0&0&0&0&0&0\\
 $\delta^{'}(\{1\},\{2\},\{3\},\{4,5\})$&1&1&1&1&0&1&1&1&0&0&1&1&0&0&0&0&0&0&0&0\\
 $\delta^{'}(\{1\},\{2\},\{3,4\},\{5\})$&1&1&1&1&0&1&1&1&0&0&0&1&0&0&0&1&0&0&0&0\\
 $\delta^{'}(\{1\},\{2\},\{3\},\{4\},\{5\})$&1&1&1&1&0&1&1&1&0&0&1&1&0&0&0&1&0&0&0&0\\
 \hline 
\end{tabular}
\caption{The representative of orbits of extreme rays in $OMCUT_5$}
\label{tab:tabl5O}
\end{table}

\begin{table}
\begin{center}
\scriptsize
\begin{tabular}{|c|c|ccccccccc|c|c|c|}
 \hline
 Or.&Representative&$O_1$&$O_2$&$O_3$&$O_4$&$O_5$&$O_6$&$O_7$&$O_8$&$O_9$&Adj.&Inc.&Size\\
 \hline
 $O_1$&$\delta^{'}(\{1\})$& 9& 20& 36& 30& 16& 54& 108& 96& 96&465&8840&10\\
 $O_2$&$\delta^{'}(\{1,5\})$& 10& 19& 38& 27& 20& 57& 90& 96& 60&417&10562&20\\
 $O_3$&$\delta^{'}(\{1\},\{2\},\{3,4,5\})$& 9& 19& 34& 24& 16& 42& 84& 72& 66&366&3106&40\\
 $O_4$&$\delta^{'}(\{1,2\},\{3\},\{4,5\})$& 10& 18& 32& 24& 16& 54& 72& 76& 40&342&3172&30\\
 $O_5$&$\delta^{'}(\{1\},\{2,3,4\},\{5\})$& 8& 20& 32& 24& 12& 42& 72& 72& 54&336&4372&20\\
 $O_6$&$\delta^{'}(\{1\},\{2,3\},\{4,5\})$& 9& 19& 28& 27& 14& 51& 64& 66& 36&314&3576&60\\
 $O_7$&$\delta^{'}(\{1\},\{2\},\{3\},\{4,5\})$& 9& 15& 28& 18& 12& 32& 57& 55& 36&262&1598&120\\
 $O_8$&$\delta^{'}(\{1\},\{2\},\{3,4\},\{5\})$& 8& 16& 24& 19& 12& 33& 55& 49& 36&252&1930&120\\
 $O_9$&$\delta^{'}(\{1\},\{2\},\{3\},\{4\},\{5\})$& 8& 10& 22& 10& 9& 18& 36& 36& 38&187&1123&120\\
 \hline
\end{tabular}
\caption{Representation matrix of $G_{OMCUT_5}$}
\label{tab:tabl5M}
\end{center}
\end{table}

There are $35320$ facets in $OMCUT_5$, which form $194$ orbits. The list of the 
representatives of these orbits is given in Appendix 2.

Besides $A_i$, $B_i$, and their zero-extensions, the simplest, i.e. with $\leq 8$ non-zeros facets of $OMCUT_5$ are 
\begin{equation*}
\begin{array}{c}
(d_{15}-d_{25})+d_{54}-d_{14}+d_{34}+d_{31}+d_{23}+d_{12}\geq 0\mbox{~and~}\\
(-d_{51}+d_{52})+d_{54}-d_{14}+d_{34}+d_{31}+d_{23}+d_{12}\geq 0\;.
\end{array}
\end{equation*}

Denote by $QHYP_n$ the cone of all quasi-semi-metrics satisfying to all hypermetric 
inequalities $H(b)$. In fact, $QHYP_n$ is polyhedral (see \cite{DGL93}) and the triangle 
inequality is redundant. The smallest case when $QHYP_n$ is a proper sub-cone of $QMET_n$ 
is $n=5$; see some information on $QHYP_5$ in Table \ref{tab:smallDatas} 
and {\it http://www.geomappl.ens.fr/\texttt{\~}dutour/}

\section{Appendix $1$: Extreme rays of the cone $QMET_5$}
By direct computation $43590$ extreme rays of $QMET_5$ were found. Under group action we got $229$ orbits.\\
See below the full list of the representatives of orbits presented as $5\times 5$ matrices. The numbers in parenthesis are, respectively, the orbit number, the adjacency of a member of the orbit, the incidence of a member of the orbit, and orbit-size. Note that adjacency is greater or equal to incidence number with equality only for the last $23$ orbits.

\newenvironment{MYmatrix}{%
\setlength{\arraycolsep}{1pt}%
\begin{matrix}}{\end{matrix}}
\newenvironment{MYarray}{%
\setlength{\arraycolsep}{1pt}%
\begin{array}}{\end{array}}
{\tiny
\begin{equation*}
\begin{MYarray}{ccccccccc}
\begin{array}{c}
(1, 8313, 64, 10)\\
\begin{MYmatrix}
0&0&0&0&1\\
0&0&0&0&1\\
0&0&0&0&1\\
0&0&0&0&1\\
0&0&0&0&0
\end{MYmatrix}
\end{array}
&\begin{array}{c}
(2, 6903, 56, 20)\\
\begin{MYmatrix}
0&0&0&1&1\\
0&0&0&1&1\\
0&0&0&1&1\\
0&0&0&0&0\\
0&0&0&0&0
\end{MYmatrix}
\end{array}
&\begin{array}{c}
(3, 2534, 48, 30)\\
\begin{MYmatrix}
0&0&1&1&1\\
0&0&1&1&1\\
0&0&0&1&1\\
0&0&0&0&0\\
0&0&0&0&0
\end{MYmatrix}
\end{array}
&\begin{array}{c}
(4, 2394, 48, 60)\\
\begin{MYmatrix}
0&0&1&1&1\\
0&0&1&1&1\\
0&0&0&0&1\\
0&0&0&0&1\\
0&0&0&0&0
\end{MYmatrix}
\end{array}
&\begin{array}{c}
(5, 1947, 52, 40)\\
\begin{MYmatrix}
0&0&0&1&1\\
0&0&0&1&1\\
0&0&0&1&1\\
0&0&0&0&1\\
0&0&0&0&0
\end{MYmatrix}
\end{array}
&\begin{array}{c}
(6, 1932, 52, 20)\\
\begin{MYmatrix}
0&1&1&1&1\\
0&0&0&0&1\\
0&0&0&0&1\\
0&0&0&0&1\\
0&0&0&0&0
\end{MYmatrix}
\end{array}
&\begin{array}{c}
(7, 1109, 44, 120)\\
\begin{MYmatrix}
0&0&1&1&1\\
0&0&1&1&1\\
0&0&0&1&1\\
0&0&0&0&1\\
0&0&0&0&0
\end{MYmatrix}
\end{array}
&\begin{array}{c}
(8, 1107, 44, 120)\\
\begin{MYmatrix}
0&1&1&1&1\\
0&0&0&1&1\\
0&0&0&1&1\\
0&0&0&0&1\\
0&0&0&0&0
\end{MYmatrix}
\end{array}
&\begin{array}{c}
(9, 701, 40, 120)\\
\begin{MYmatrix}
0&1&1&1&1\\
0&0&1&1&1\\
0&0&0&1&1\\
0&0&0&0&1\\
0&0&0&0&0
\end{MYmatrix}
\end{array}
\end{MYarray}
\end{equation*}
\begin{equation*}
\begin{MYarray}{ccccccccc}
\begin{array}{c}
(10, 461, 38, 120)\\
\begin{MYmatrix}
0&1&1&1&1\\
0&0&1&1&1\\
0&0&0&0&1\\
0&0&0&0&1\\
0&0&1&1&0
\end{MYmatrix}
\end{array}
&\begin{array}{c}
(11, 444, 36, 120)\\
\begin{MYmatrix}
0&1&1&1&1\\
0&0&1&1&1\\
0&0&0&1&1\\
0&0&0&0&1\\
0&0&0&1&0
\end{MYmatrix}
\end{array}
&\begin{array}{c}
(12, 413, 40, 60)\\
\begin{MYmatrix}
0&1&1&1&1\\
0&0&0&1&1\\
0&0&0&1&1\\
0&0&0&0&1\\
0&0&0&1&0
\end{MYmatrix}
\end{array}
&\begin{array}{c}
(13, 404, 41, 120)\\
\begin{MYmatrix}
0&1&1&1&2\\
0&0&0&1&1\\
0&0&0&1&1\\
0&0&0&0&1\\
0&0&0&0&0
\end{MYmatrix}
\end{array}
&\begin{array}{c}
(14, 401, 40, 60)\\
\begin{MYmatrix}
0&0&1&1&1\\
0&0&1&1&1\\
0&0&0&1&1\\
0&0&0&0&1\\
0&0&0&1&0
\end{MYmatrix}
\end{array}
&\begin{array}{c}
(15, 398, 40, 120)\\
\begin{MYmatrix}
0&0&1&1&2\\
0&0&1&1&2\\
0&0&0&1&1\\
0&0&0&0&1\\
0&0&0&0&0
\end{MYmatrix}
\end{array}
&\begin{array}{c}
(16, 376, 32, 30)\\
\begin{MYmatrix}
0&1&1&1&1\\
0&0&1&1&0\\
0&0&0&1&0\\
0&0&1&0&0\\
1&1&1&1&0
\end{MYmatrix}
\end{array}
&\begin{array}{c}
(17, 355, 36, 120)\\
\begin{MYmatrix}
0&1&1&1&1\\
0&0&1&1&1\\
0&0&0&1&1\\
0&0&0&0&0\\
0&0&1&1&0
\end{MYmatrix}
\end{array}
&\begin{array}{c}
(18, 345, 37, 120)\\
\begin{MYmatrix}
0&1&1&1&2\\
0&0&1&1&1\\
0&0&0&1&1\\
0&0&0&0&1\\
0&0&0&0&0
\end{MYmatrix}
\end{array}
\end{MYarray}
\end{equation*}
\begin{equation*}
\begin{MYarray}{ccccccccc}
\begin{array}{c}
(19, 342, 40, 60)\\
\begin{MYmatrix}
0&0&1&1&1\\
0&0&1&1&1\\
0&0&0&1&1\\
0&0&0&0&0\\
0&0&1&1&0
\end{MYmatrix}
\end{array}
&\begin{array}{c}
(20, 316, 38, 60)\\
\begin{MYmatrix}
0&1&1&1&1\\
0&0&0&1&1\\
0&0&0&1&1\\
0&0&0&0&0\\
0&1&1&1&0
\end{MYmatrix}
\end{array}
&\begin{array}{c}
(21, 292, 34, 60)\\
\begin{MYmatrix}
0&0&1&1&1\\
0&0&1&1&1\\
0&0&0&1&0\\
0&0&1&0&0\\
1&1&1&1&0
\end{MYmatrix}
\end{array}
&\begin{array}{c}
(22, 287, 33, 120)\\
\begin{MYmatrix}
0&1&1&2&1\\
0&0&1&1&1\\
0&0&0&1&1\\
0&0&0&0&0\\
0&0&1&1&0
\end{MYmatrix}
\end{array}
&\begin{array}{c}
(23, 286, 32, 60)\\
\begin{MYmatrix}
0&1&1&1&1\\
0&0&1&1&1\\
0&0&0&1&0\\
0&0&1&0&0\\
0&1&1&1&0
\end{MYmatrix}
\end{array}
&\begin{array}{c}
(24, 248, 29, 120)\\
\begin{MYmatrix}
0&1&1&1&1\\
0&0&1&1&1\\
0&0&0&1&0\\
0&0&1&0&0\\
1&1&1&1&0
\end{MYmatrix}
\end{array}
&\begin{array}{c}
(25, 186, 33, 240)\\
\begin{MYmatrix}
0&1&1&1&1\\
0&0&1&1&1\\
0&0&0&1&1\\
0&0&0&0&1\\
0&0&1&1&0
\end{MYmatrix}
\end{array}
&\begin{array}{c}
(26, 175, 34, 240)\\
\begin{MYmatrix}
0&1&1&1&2\\
0&0&1&1&1\\
0&0&0&1&1\\
0&0&0&0&1\\
0&0&0&1&0
\end{MYmatrix}
\end{array}
&\begin{array}{c}
(27, 175, 30, 120)\\
\begin{MYmatrix}
0&1&1&2&1\\
0&0&1&1&1\\
0&0&0&1&0\\
0&0&1&0&0\\
0&1&1&1&0
\end{MYmatrix}
\end{array}
\end{MYarray}
\end{equation*}
\begin{equation*}
\begin{MYarray}{ccccccccc}
\begin{array}{c}
(28, 174, 31, 120)\\
\begin{MYmatrix}
0&1&1&2&1\\
0&0&1&2&1\\
0&0&0&1&1\\
0&0&0&0&0\\
0&0&1&1&0
\end{MYmatrix}
\end{array}
&\begin{array}{c}
(29, 172, 34, 120)\\
\begin{MYmatrix}
0&0&1&1&1\\
0&0&1&1&1\\
0&0&0&1&0\\
0&0&1&0&0\\
1&1&1&2&0
\end{MYmatrix}
\end{array}
&\begin{array}{c}
(30, 172, 34, 120)\\
\begin{MYmatrix}
0&0&1&2&1\\
0&0&1&2&1\\
0&0&0&1&0\\
0&0&1&0&0\\
1&1&1&1&0
\end{MYmatrix}
\end{array}
&\begin{array}{c}
(31, 160, 30, 120)\\
\begin{MYmatrix}
0&1&1&2&1\\
0&0&1&1&1\\
0&0&0&1&1\\
0&0&0&0&0\\
0&1&1&1&0
\end{MYmatrix}
\end{array}
&\begin{array}{c}
(32, 149, 35, 240)\\
\begin{MYmatrix}
0&1&1&1&2\\
0&0&1&1&2\\
0&0&0&1&1\\
0&0&0&0&1\\
0&0&0&0&0
\end{MYmatrix}
\end{array}
&\begin{array}{c}
(33, 143, 28, 240)\\
\begin{MYmatrix}
0&1&1&2&1\\
0&0&1&1&1\\
0&0&0&1&0\\
0&0&1&0&0\\
1&1&1&1&0
\end{MYmatrix}
\end{array}
&\begin{array}{c}
(34, 138, 29, 240)\\
\begin{MYmatrix}
0&1&1&1&1\\
0&0&1&1&1\\
0&0&0&1&0\\
0&0&1&0&0\\
1&1&1&2&0
\end{MYmatrix}
\end{array}
&\begin{array}{c}
(35, 138, 24, 20)\\
\begin{MYmatrix}
0&1&1&1&1\\
0&0&1&1&0\\
0&1&0&1&0\\
0&1&1&0&0\\
1&1&1&1&0
\end{MYmatrix}
\end{array}
&\begin{array}{c}
(36, 133, 26, 60)\\
\begin{MYmatrix}
0&1&2&2&1\\
0&0&1&1&1\\
0&0&0&1&0\\
0&0&1&0&0\\
0&1&1&1&0
\end{MYmatrix}
\end{array}
\end{MYarray}
\end{equation*}
\begin{equation*}
\begin{MYarray}{ccccccccc}
\begin{array}{c}
(37, 131, 25, 20)\\
\begin{MYmatrix}
0&1&2&1&1\\
0&0&1&1&1\\
0&0&0&0&0\\
0&1&1&0&1\\
0&1&1&1&0
\end{MYmatrix}
\end{array}
&\begin{array}{c}
(38, 128, 33, 120)\\
\begin{MYmatrix}
0&1&1&1&1\\
0&0&1&1&1\\
0&0&0&1&1\\
0&0&0&0&0\\
0&1&1&1&0
\end{MYmatrix}
\end{array}
&\begin{array}{c}
(39, 122, 26, 120)\\
\begin{MYmatrix}
0&1&2&1&1\\
0&0&1&1&1\\
0&0&0&1&0\\
0&1&1&0&0\\
1&1&1&1&0
\end{MYmatrix}
\end{array}
&\begin{array}{c}
(40, 116, 27, 120)\\
\begin{MYmatrix}
0&1&1&1&1\\
0&0&1&1&1\\
0&0&0&1&0\\
0&1&1&0&0\\
1&1&1&1&0
\end{MYmatrix}
\end{array}
&\begin{array}{c}
(41, 113, 24, 240)\\
\begin{MYmatrix}
0&1&2&2&1\\
0&0&1&1&1\\
0&0&0&1&0\\
0&1&1&0&0\\
1&1&1&1&0
\end{MYmatrix}
\end{array}
&\begin{array}{c}
(42, 113, 25, 120)\\
\begin{MYmatrix}
0&1&2&2&1\\
0&0&1&1&1\\
0&0&0&1&0\\
0&0&1&0&0\\
1&1&1&1&0
\end{MYmatrix}
\end{array}
&\begin{array}{c}
(43, 109, 27, 240)\\
\begin{MYmatrix}
0&1&1&2&1\\
0&0&1&1&1\\
0&0&0&1&0\\
0&1&1&0&0\\
1&1&1&1&0
\end{MYmatrix}
\end{array}
&\begin{array}{c}
(44, 107, 31, 120)\\
\begin{MYmatrix}
0&1&1&1&1\\
0&0&1&1&0\\
0&0&0&1&0\\
0&0&1&0&0\\
1&1&1&2&0
\end{MYmatrix}
\end{array}
&\begin{array}{c}
(45, 102, 25, 240)\\
\begin{MYmatrix}
0&1&2&2&2\\
0&0&2&1&1\\
0&0&0&1&0\\
0&1&1&0&0\\
1&1&1&1&0
\end{MYmatrix}
\end{array}
\end{MYarray}
\end{equation*}
\begin{equation*}
\begin{MYarray}{ccccccccc}
\begin{array}{c}
(46, 98, 29, 120)\\
\begin{MYmatrix}
0&1&2&2&1\\
0&0&1&1&0\\
0&1&0&2&1\\
0&1&0&0&0\\
1&1&1&1&0
\end{MYmatrix}
\end{array}
&\begin{array}{c}
(47, 95, 31, 240)\\
\begin{MYmatrix}
0&1&1&2&1\\
0&0&1&2&1\\
0&0&0&2&1\\
0&0&0&0&1\\
0&0&1&1&0
\end{MYmatrix}
\end{array}
&\begin{array}{c}
(48, 91, 30, 240)\\
\begin{MYmatrix}
0&1&1&2&1\\
0&0&1&2&1\\
0&0&0&1&1\\
0&0&0&0&1\\
0&0&1&1&0
\end{MYmatrix}
\end{array}
&\begin{array}{c}
(49, 91, 24, 120)\\
\begin{MYmatrix}
0&1&1&1&1\\
0&0&1&1&0\\
0&1&0&1&0\\
0&1&1&0&0\\
1&1&1&2&0
\end{MYmatrix}
\end{array}
&\begin{array}{c}
(50, 91, 30, 120)\\
\begin{MYmatrix}
0&1&1&2&2\\
0&0&1&1&1\\
0&0&0&1&1\\
0&0&0&0&1\\
0&0&0&1&0
\end{MYmatrix}
\end{array}
&\begin{array}{c}
(51, 90, 26, 120)\\
\begin{MYmatrix}
0&1&2&1&1\\
0&0&2&1&0\\
0&0&0&1&0\\
0&1&1&0&0\\
1&1&2&1&0
\end{MYmatrix}
\end{array}
&\begin{array}{c}
(52, 88, 28, 40)\\
\begin{MYmatrix}
0&1&1&1&1\\
0&0&1&1&1\\
0&0&0&1&1\\
0&0&1&0&1\\
0&0&1&1&0
\end{MYmatrix}
\end{array}
&\begin{array}{c}
(53, 87, 31, 240)\\
\begin{MYmatrix}
0&1&1&2&1\\
0&0&1&1&1\\
0&0&0&1&1\\
0&0&0&0&1\\
0&0&1&1&0
\end{MYmatrix}
\end{array}
&\begin{array}{c}
(54, 84, 21, 40)\\
\begin{MYmatrix}
0&2&2&2&2\\
0&0&2&1&0\\
0&1&0&2&0\\
0&2&1&0&0\\
2&2&2&2&0
\end{MYmatrix}
\end{array}
\end{MYarray}
\end{equation*}
\begin{equation*}
\begin{MYarray}{ccccccccc}
\begin{array}{c}
(55, 83, 32, 240)\\
\begin{MYmatrix}
0&1&1&2&2\\
0&0&1&1&2\\
0&0&0&1&2\\
0&0&0&0&1\\
0&0&0&0&0
\end{MYmatrix}
\end{array}
&\begin{array}{c}
(56, 83, 33, 240)\\
\begin{MYmatrix}
0&1&1&2&2\\
0&0&1&1&2\\
0&0&0&1&2\\
0&0&0&0&2\\
0&0&0&0&0
\end{MYmatrix}
\end{array}
&\begin{array}{c}
(57, 81, 31, 240)\\
\begin{MYmatrix}
0&1&1&2&1\\
0&0&1&2&1\\
0&0&0&1&1\\
0&0&0&0&0\\
0&0&1&2&0
\end{MYmatrix}
\end{array}
&\begin{array}{c}
(58, 75, 30, 240)\\
\begin{MYmatrix}
0&1&1&2&1\\
0&0&1&1&1\\
0&0&0&1&0\\
0&0&1&0&0\\
0&1&1&2&0
\end{MYmatrix}
\end{array}
&\begin{array}{c}
(59, 75, 27, 120)\\
\begin{MYmatrix}
0&1&1&2&1\\
0&0&1&1&1\\
0&0&0&1&0\\
0&1&1&0&0\\
1&1&2&1&0
\end{MYmatrix}
\end{array}
&\begin{array}{c}
(60, 72, 27, 120)\\
\begin{MYmatrix}
0&1&2&2&2\\
0&0&1&2&1\\
0&0&0&1&1\\
0&0&0&0&0\\
0&0&1&1&0
\end{MYmatrix}
\end{array}
&\begin{array}{c}
(61, 71, 25, 240)\\
\begin{MYmatrix}
0&1&2&2&1\\
0&0&2&1&1\\
0&0&0&1&0\\
0&1&1&0&0\\
1&1&1&1&0
\end{MYmatrix}
\end{array}
&\begin{array}{c}
(62, 70, 25, 240)\\
\begin{MYmatrix}
0&1&2&2&1\\
0&0&1&1&1\\
0&0&0&1&0\\
0&0&1&0&1\\
0&1&1&1&0
\end{MYmatrix}
\end{array}
&\begin{array}{c}
(63, 70, 22, 240)\\
\begin{MYmatrix}
0&1&2&2&2\\
0&0&2&1&1\\
0&0&0&1&0\\
0&1&1&0&0\\
1&1&2&1&0
\end{MYmatrix}
\end{array}
\end{MYarray}
\end{equation*}
\begin{equation*}
\begin{MYarray}{ccccccccc}
\begin{array}{c}
(64, 69, 26, 240)\\
\begin{MYmatrix}
0&1&2&2&2\\
0&0&1&2&1\\
0&0&0&1&1\\
0&0&0&0&0\\
0&1&1&1&0
\end{MYmatrix}
\end{array}
&\begin{array}{c}
(65, 68, 28, 240)\\
\begin{MYmatrix}
0&1&1&1&1\\
0&0&1&1&1\\
0&0&0&1&0\\
0&0&1&0&1\\
1&1&1&1&0
\end{MYmatrix}
\end{array}
&\begin{array}{c}
(66, 68, 27, 240)\\
\begin{MYmatrix}
0&1&2&1&1\\
0&0&2&1&1\\
0&0&0&1&0\\
0&1&1&0&0\\
1&1&1&1&0
\end{MYmatrix}
\end{array}
&\begin{array}{c}
(67, 65, 30, 240)\\
\begin{MYmatrix}
0&1&1&2&1\\
0&0&1&2&1\\
0&0&0&1&0\\
0&0&1&0&1\\
0&1&1&1&0
\end{MYmatrix}
\end{array}
&\begin{array}{c}
(68, 64, 26, 240)\\
\begin{MYmatrix}
0&1&2&2&1\\
0&0&1&1&1\\
0&0&0&1&0\\
0&0&1&0&0\\
0&1&1&2&0
\end{MYmatrix}
\end{array}
&\begin{array}{c}
(69, 64, 29, 240)\\
\begin{MYmatrix}
0&1&1&1&1\\
0&0&1&1&1\\
0&0&0&1&1\\
0&0&1&0&0\\
1&1&1&2&0
\end{MYmatrix}
\end{array}
&\begin{array}{c}
(70, 64, 26, 240)\\
\begin{MYmatrix}
0&1&2&2&2\\
0&0&1&2&1\\
0&0&0&1&0\\
0&0&1&0&0\\
0&1&1&1&0
\end{MYmatrix}
\end{array}
&\begin{array}{c}
(71, 64, 29, 240)\\
\begin{MYmatrix}
0&1&1&2&1\\
0&0&1&1&1\\
0&0&0&1&0\\
0&0&1&0&1\\
0&1&1&1&0
\end{MYmatrix}
\end{array}
&\begin{array}{c}
(72, 64, 29, 240)\\
\begin{MYmatrix}
0&1&1&2&1\\
0&0&1&1&1\\
0&0&0&1&1\\
0&0&0&0&1\\
0&1&1&1&0
\end{MYmatrix}
\end{array}
\end{MYarray}
\end{equation*}
\begin{equation*}
\begin{MYarray}{ccccccccc}
\begin{array}{c}
(73, 64, 30, 240)\\
\begin{MYmatrix}
0&1&1&2&2\\
0&0&1&1&1\\
0&0&0&1&1\\
0&0&0&0&0\\
0&1&1&1&0
\end{MYmatrix}
\end{array}
&\begin{array}{c}
(74, 63, 24, 120)\\
\begin{MYmatrix}
0&1&2&1&1\\
0&0&1&1&1\\
0&0&0&0&1\\
0&1&1&0&1\\
0&1&1&1&0
\end{MYmatrix}
\end{array}
&\begin{array}{c}
(75, 63, 24, 120)\\
\begin{MYmatrix}
0&1&1&2&1\\
0&0&1&1&0\\
0&1&0&1&0\\
0&1&1&0&1\\
1&1&1&1&0
\end{MYmatrix}
\end{array}
&\begin{array}{c}
(76, 63, 28, 120)\\
\begin{MYmatrix}
0&1&1&1&1\\
0&0&1&1&1\\
0&0&0&1&0\\
0&1&1&0&0\\
1&1&1&2&0
\end{MYmatrix}
\end{array}
&\begin{array}{c}
(77, 63, 28, 120)\\
\begin{MYmatrix}
0&1&2&2&1\\
0&0&1&2&1\\
0&0&0&1&1\\
0&0&0&0&0\\
0&1&1&1&0
\end{MYmatrix}
\end{array}
&\begin{array}{c}
(78, 62, 28, 240)\\
\begin{MYmatrix}
0&1&1&2&1\\
0&0&1&2&1\\
0&0&0&1&0\\
0&0&1&0&0\\
1&1&1&1&0
\end{MYmatrix}
\end{array}
&\begin{array}{c}
(79, 62, 25, 240)\\
\begin{MYmatrix}
0&1&2&2&1\\
0&0&1&2&1\\
0&0&0&1&1\\
0&0&1&0&0\\
0&1&1&1&0
\end{MYmatrix}
\end{array}
&\begin{array}{c}
(80, 61, 28, 240)\\
\begin{MYmatrix}
0&1&1&2&1\\
0&0&1&1&1\\
0&0&0&1&0\\
0&0&1&0&0\\
1&1&2&1&0
\end{MYmatrix}
\end{array}
&\begin{array}{c}
(81, 60, 32, 240)\\
\begin{MYmatrix}
0&1&1&2&2\\
0&0&1&1&1\\
0&0&0&1&1\\
0&0&0&0&0\\
0&0&1&1&0
\end{MYmatrix}
\end{array}
\end{MYarray}
\end{equation*}
\begin{equation*}
\begin{MYarray}{ccccccccc}
\begin{array}{c}
(82, 58, 30, 240)\\
\begin{MYmatrix}
0&1&1&1&1\\
0&0&1&1&1\\
0&0&0&1&0\\
0&0&1&0&1\\
0&1&1&1&0
\end{MYmatrix}
\end{array}
&\begin{array}{c}
(83, 58, 23, 240)\\
\begin{MYmatrix}
0&1&1&2&1\\
0&0&1&1&0\\
0&1&0&1&1\\
0&1&1&0&0\\
1&1&1&1&0
\end{MYmatrix}
\end{array}
&\begin{array}{c}
(84, 57, 33, 120)\\
\begin{MYmatrix}
0&1&1&2&2\\
0&0&1&1&2\\
0&0&0&1&1\\
0&0&0&0&1\\
0&0&0&0&0
\end{MYmatrix}
\end{array}
&\begin{array}{c}
(85, 54, 28, 240)\\
\begin{MYmatrix}
0&1&2&2&2\\
0&0&1&2&1\\
0&0&0&1&1\\
0&0&0&0&0\\
0&0&1&2&0
\end{MYmatrix}
\end{array}
&\begin{array}{c}
(86, 54, 25, 120)\\
\begin{MYmatrix}
0&1&2&1&2\\
0&0&1&1&1\\
0&0&0&0&0\\
0&1&1&0&1\\
0&1&1&1&0
\end{MYmatrix}
\end{array}
&\begin{array}{c}
(87, 54, 25, 120)\\
\begin{MYmatrix}
0&1&2&2&1\\
0&0&1&1&1\\
0&0&0&0&0\\
0&1&1&0&1\\
0&1&2&1&0
\end{MYmatrix}
\end{array}
&\begin{array}{c}
(88, 53, 24, 240)\\
\begin{MYmatrix}
0&1&2&2&1\\
0&0&1&2&1\\
0&0&0&1&1\\
0&0&1&0&0\\
0&1&1&2&0
\end{MYmatrix}
\end{array}
&\begin{array}{c}
(89, 53, 27, 240)\\
\begin{MYmatrix}
0&1&2&1&1\\
0&0&2&1&0\\
0&0&0&1&0\\
0&1&1&0&0\\
1&2&2&1&0
\end{MYmatrix}
\end{array}
&\begin{array}{c}
(90, 51, 25, 240)\\
\begin{MYmatrix}
0&1&2&2&2\\
0&0&1&2&1\\
0&0&0&1&0\\
0&0&1&0&0\\
0&1&1&2&0
\end{MYmatrix}
\end{array}
\end{MYarray}
\end{equation*}
\begin{equation*}
\begin{MYarray}{ccccccccc}
\begin{array}{c}
(91, 51, 27, 240)\\
\begin{MYmatrix}
0&1&2&2&2\\
0&0&1&2&2\\
0&0&0&2&1\\
0&0&0&0&0\\
0&1&1&1&0
\end{MYmatrix}
\end{array}
&\begin{array}{c}
(92, 51, 25, 240)\\
\begin{MYmatrix}
0&1&2&1&1\\
0&0&2&1&0\\
0&0&0&1&0\\
0&1&2&0&0\\
1&2&2&1&0
\end{MYmatrix}
\end{array}
&\begin{array}{c}
(93, 51, 28, 20)\\
\begin{MYmatrix}
0&1&1&1&1\\
0&0&1&1&1\\
0&0&0&0&0\\
0&1&1&0&1\\
0&1&1&1&0
\end{MYmatrix}
\end{array}
&\begin{array}{c}
(94, 48, 25, 80)\\
\begin{MYmatrix}
0&2&2&2&2\\
0&0&2&2&2\\
0&0&0&2&1\\
0&0&1&0&2\\
0&0&2&1&0
\end{MYmatrix}
\end{array}
&\begin{array}{c}
(95, 47, 26, 240)\\
\begin{MYmatrix}
0&1&2&2&1\\
0&0&1&1&1\\
0&0&0&1&1\\
0&0&1&0&0\\
0&1&1&2&0
\end{MYmatrix}
\end{array}
&\begin{array}{c}
(96, 47, 26, 240)\\
\begin{MYmatrix}
0&1&1&2&1\\
0&0&1&1&1\\
0&0&0&1&0\\
0&0&1&0&0\\
1&1&1&2&0
\end{MYmatrix}
\end{array}
&\begin{array}{c}
(97, 47, 26, 240)\\
\begin{MYmatrix}
0&1&2&2&2\\
0&0&1&2&1\\
0&0&0&1&0\\
0&1&1&0&0\\
0&1&2&1&0
\end{MYmatrix}
\end{array}
&\begin{array}{c}
(98, 47, 27, 240)\\
\begin{MYmatrix}
0&1&2&3&2\\
0&0&2&2&1\\
0&0&0&1&1\\
0&0&0&0&0\\
0&1&1&1&0
\end{MYmatrix}
\end{array}
&\begin{array}{c}
(99, 47, 29, 120)\\
\begin{MYmatrix}
0&1&2&2&3\\
0&0&1&2&2\\
0&0&0&1&2\\
0&0&0&0&1\\
0&0&0&0&0
\end{MYmatrix}
\end{array}
\end{MYarray}
\end{equation*}
\begin{equation*}
\begin{MYarray}{ccccccccc}
\begin{array}{c}
(100, 47, 26, 120)\\
\begin{MYmatrix}
0&1&2&2&2\\
0&0&1&1&1\\
0&0&0&1&0\\
0&0&1&0&0\\
0&1&1&1&0
\end{MYmatrix}
\end{array}
&\begin{array}{c}
(101, 47, 23, 120)\\
\begin{MYmatrix}
0&1&1&1&1\\
0&0&1&1&0\\
0&1&0&1&0\\
0&1&1&0&1\\
1&1&1&1&0
\end{MYmatrix}
\end{array}
&\begin{array}{c}
(102, 47, 31, 120)\\
\begin{MYmatrix}
0&1&1&2&2\\
0&0&1&1&1\\
0&0&0&1&1\\
0&0&0&0&0\\
0&1&1&2&0
\end{MYmatrix}
\end{array}
&\begin{array}{c}
(103, 46, 25, 240)\\
\begin{MYmatrix}
0&2&2&2&2\\
0&0&1&1&1\\
0&0&0&1&0\\
0&0&1&0&0\\
0&1&1&2&0
\end{MYmatrix}
\end{array}
&\begin{array}{c}
(104, 45, 22, 240)\\
\begin{MYmatrix}
0&1&1&2&1\\
0&0&1&1&0\\
0&1&0&1&1\\
0&1&1&0&0\\
1&1&1&2&0
\end{MYmatrix}
\end{array}
&\begin{array}{c}
(105, 45, 22, 120)\\
\begin{MYmatrix}
0&1&2&2&1\\
0&0&1&1&1\\
0&0&0&1&0\\
0&1&1&0&0\\
1&1&2&1&0
\end{MYmatrix}
\end{array}
&\begin{array}{c}
(106, 44, 26, 240)\\
\begin{MYmatrix}
0&1&1&2&1\\
0&0&1&1&1\\
0&0&0&1&0\\
0&1&1&0&0\\
1&1&1&2&0
\end{MYmatrix}
\end{array}
&\begin{array}{c}
(107, 44, 20, 240)\\
\begin{MYmatrix}
0&2&2&2&2\\
0&0&2&1&0\\
0&1&0&2&0\\
0&2&1&0&0\\
2&2&2&3&0
\end{MYmatrix}
\end{array}
&\begin{array}{c}
(108, 44, 20, 120)\\
\begin{MYmatrix}
0&2&3&2&1\\
0&0&1&2&0\\
0&1&0&1&1\\
0&2&1&0&0\\
2&2&2&2&0
\end{MYmatrix}
\end{array}
\end{MYarray}
\end{equation*}
\begin{equation*}
\begin{MYarray}{ccccccccc}
\begin{array}{c}
(109, 43, 25, 240)\\
\begin{MYmatrix}
0&1&2&2&1\\
0&0&1&1&0\\
0&1&0&2&1\\
0&1&0&0&0\\
1&1&1&2&0
\end{MYmatrix}
\end{array}
&\begin{array}{c}
(110, 43, 25, 240)\\
\begin{MYmatrix}
0&1&2&2&1\\
0&0&1&1&1\\
0&0&0&1&0\\
0&0&1&0&1\\
1&1&1&1&0
\end{MYmatrix}
\end{array}
&\begin{array}{c}
(111, 42, 26, 240)\\
\begin{MYmatrix}
0&2&2&2&1\\
0&0&1&2&1\\
0&0&0&2&0\\
0&0&1&0&0\\
1&1&1&3&0
\end{MYmatrix}
\end{array}
&\begin{array}{c}
(112, 42, 23, 120)\\
\begin{MYmatrix}
0&1&2&1&1\\
0&0&1&1&1\\
0&0&0&1&0\\
0&1&1&0&1\\
1&1&1&1&0
\end{MYmatrix}
\end{array}
&\begin{array}{c}
(113, 41, 27, 240)\\
\begin{MYmatrix}
0&1&1&2&1\\
0&0&1&2&1\\
0&0&0&1&0\\
0&1&1&0&0\\
1&1&2&1&0
\end{MYmatrix}
\end{array}
&\begin{array}{c}
(114, 40, 26, 240)\\
\begin{MYmatrix}
0&1&3&2&2\\
0&0&2&1&1\\
0&0&0&1&0\\
0&1&1&0&0\\
1&1&1&1&0
\end{MYmatrix}
\end{array}
&\begin{array}{c}
(115, 40, 26, 240)\\
\begin{MYmatrix}
0&1&2&2&1\\
0&0&1&2&1\\
0&0&0&2&1\\
0&0&1&0&0\\
0&1&1&1&0
\end{MYmatrix}
\end{array}
&\begin{array}{c}
(116, 40, 26, 240)\\
\begin{MYmatrix}
0&1&2&1&1\\
0&0&2&1&0\\
0&0&0&1&0\\
0&1&1&0&0\\
1&1&3&2&0
\end{MYmatrix}
\end{array}
&\begin{array}{c}
(117, 39, 23, 240)\\
\begin{MYmatrix}
0&1&2&2&2\\
0&0&1&2&1\\
0&0&0&1&0\\
0&1&1&0&0\\
1&1&2&1&0
\end{MYmatrix}
\end{array}
\end{MYarray}
\end{equation*}
\begin{equation*}
\begin{MYarray}{ccccccccc}
\begin{array}{c}
(118, 39, 26, 240)\\
\begin{MYmatrix}
0&1&2&2&2\\
0&0&1&2&1\\
0&0&0&2&0\\
0&0&1&0&0\\
0&1&1&2&0
\end{MYmatrix}
\end{array}
&\begin{array}{c}
(119, 39, 23, 240)\\
\begin{MYmatrix}
0&1&3&2&1\\
0&0&2&1&1\\
0&0&0&1&0\\
0&1&1&0&0\\
1&1&2&1&0
\end{MYmatrix}
\end{array}
&\begin{array}{c}
(120, 39, 23, 240)\\
\begin{MYmatrix}
0&1&1&2&1\\
0&0&1&1&1\\
0&1&0&2&0\\
0&1&1&0&0\\
1&1&1&2&0
\end{MYmatrix}
\end{array}
&\begin{array}{c}
(121, 39, 23, 240)\\
\begin{MYmatrix}
0&1&2&2&2\\
0&0&1&1&1\\
0&0&0&1&0\\
0&1&1&0&0\\
1&1&2&1&0
\end{MYmatrix}
\end{array}
&\begin{array}{c}
(122, 38, 24, 240)\\
\begin{MYmatrix}
0&2&2&2&1\\
0&0&1&2&1\\
0&0&0&2&0\\
0&0&1&0&0\\
1&1&1&2&0
\end{MYmatrix}
\end{array}
&\begin{array}{c}
(123, 38, 27, 240)\\
\begin{MYmatrix}
0&1&2&2&2\\
0&0&1&1&1\\
0&0&0&1&0\\
0&0&1&0&0\\
0&1&1&2&0
\end{MYmatrix}
\end{array}
&\begin{array}{c}
(124, 38, 24, 240)\\
\begin{MYmatrix}
0&1&1&1&1\\
0&0&1&1&0\\
0&1&0&1&1\\
0&1&1&0&0\\
1&1&1&2&0
\end{MYmatrix}
\end{array}
&\begin{array}{c}
(125, 38, 27, 240)\\
\begin{MYmatrix}
0&1&2&3&1\\
0&0&1&2&1\\
0&0&0&2&0\\
0&0&0&0&0\\
0&1&2&2&0
\end{MYmatrix}
\end{array}
&\begin{array}{c}
(126, 38, 24, 240)\\
\begin{MYmatrix}
0&1&2&2&1\\
0&0&1&2&1\\
0&0&0&1&1\\
0&1&1&0&0\\
1&1&1&1&0
\end{MYmatrix}
\end{array}
\end{MYarray}
\end{equation*}
\begin{equation*}
\begin{MYarray}{ccccccccc}
\begin{array}{c}
(127, 38, 25, 120)\\
\begin{MYmatrix}
0&1&2&2&2\\
0&0&1&2&1\\
0&0&0&1&1\\
0&0&0&0&0\\
0&1&1&2&0
\end{MYmatrix}
\end{array}
&\begin{array}{c}
(128, 37, 24, 240)\\
\begin{MYmatrix}
0&1&2&1&2\\
0&0&2&1&1\\
0&0&0&1&0\\
0&1&1&0&1\\
1&1&1&1&0
\end{MYmatrix}
\end{array}
&\begin{array}{c}
(129, 37, 26, 240)\\
\begin{MYmatrix}
0&1&2&2&1\\
0&0&1&2&1\\
0&0&0&1&1\\
0&0&0&0&1\\
0&1&1&2&0
\end{MYmatrix}
\end{array}
&\begin{array}{c}
(130, 37, 24, 240)\\
\begin{MYmatrix}
0&1&2&1&1\\
0&0&2&1&0\\
0&0&0&1&0\\
0&1&1&0&0\\
1&1&2&2&0
\end{MYmatrix}
\end{array}
&\begin{array}{c}
(131, 36, 25, 240)\\
\begin{MYmatrix}
0&1&1&2&1\\
0&0&1&1&0\\
0&1&0&1&0\\
0&1&1&0&1\\
1&1&2&1&0
\end{MYmatrix}
\end{array}
&\begin{array}{c}
(132, 36, 25, 240)\\
\begin{MYmatrix}
0&1&2&1&1\\
0&0&1&1&1\\
0&0&0&1&0\\
0&1&1&0&1\\
0&1&2&1&0
\end{MYmatrix}
\end{array}
&\begin{array}{c}
(133, 35, 25, 240)\\
\begin{MYmatrix}
0&2&2&2&2\\
0&0&1&2&2\\
0&0&0&2&1\\
0&0&1&0&0\\
0&1&2&1&0
\end{MYmatrix}
\end{array}
&\begin{array}{c}
(134, 35, 25, 240)\\
\begin{MYmatrix}
0&1&2&2&2\\
0&0&1&1&2\\
0&0&0&0&1\\
0&1&1&0&2\\
0&1&1&0&0
\end{MYmatrix}
\end{array}
&\begin{array}{c}
(135, 35, 27, 240)\\
\begin{MYmatrix}
0&1&1&2&1\\
0&0&1&2&1\\
0&0&0&1&0\\
0&1&1&0&0\\
1&1&1&1&0
\end{MYmatrix}
\end{array}
\end{MYarray}
\end{equation*}
\begin{equation*}
\begin{MYarray}{ccccccccc}
\begin{array}{c}
(136, 35, 26, 240)\\
\begin{MYmatrix}
0&1&2&3&3\\
0&0&2&2&2\\
0&0&0&2&1\\
0&0&0&0&0\\
0&1&1&1&0
\end{MYmatrix}
\end{array}
&\begin{array}{c}
(137, 34, 24, 240)\\
\begin{MYmatrix}
0&1&1&2&1\\
0&0&1&1&0\\
0&1&0&2&1\\
0&1&1&0&0\\
1&1&1&1&0
\end{MYmatrix}
\end{array}
&\begin{array}{c}
(138, 34, 27, 240)\\
\begin{MYmatrix}
0&1&2&3&2\\
0&0&1&2&1\\
0&0&0&1&0\\
0&0&1&0&0\\
0&1&1&1&0
\end{MYmatrix}
\end{array}
&\begin{array}{c}
(139, 34, 24, 120)\\
\begin{MYmatrix}
0&1&1&2&1\\
0&0&1&1&0\\
0&1&0&1&0\\
0&1&1&0&0\\
1&1&2&1&0
\end{MYmatrix}
\end{array}
&\begin{array}{c}
(140, 32, 26, 240)\\
\begin{MYmatrix}
0&1&2&2&2\\
0&0&1&2&2\\
0&0&0&1&1\\
0&0&0&0&0\\
0&1&1&2&0
\end{MYmatrix}
\end{array}
&\begin{array}{c}
(141, 32, 25, 240)\\
\begin{MYmatrix}
0&1&2&2&1\\
0&0&1&2&0\\
0&1&0&2&1\\
0&1&0&0&0\\
1&1&1&3&0
\end{MYmatrix}
\end{array}
&\begin{array}{c}
(142, 32, 23, 240)\\
\begin{MYmatrix}
0&1&2&2&1\\
0&0&1&1&1\\
0&0&0&1&0\\
0&1&1&0&1\\
0&1&2&1&0
\end{MYmatrix}
\end{array}
&\begin{array}{c}
(143, 32, 25, 240)\\
\begin{MYmatrix}
0&1&2&3&1\\
0&0&1&2&1\\
0&0&0&1&1\\
0&0&1&0&0\\
0&1&1&2&0
\end{MYmatrix}
\end{array}
&\begin{array}{c}
(144, 32, 26, 240)\\
\begin{MYmatrix}
0&2&2&3&2\\
0&0&1&1&1\\
0&0&0&1&0\\
0&0&1&0&0\\
0&1&2&1&0
\end{MYmatrix}
\end{array}
\end{MYarray}
\end{equation*}
\begin{equation*}
\begin{MYarray}{ccccccccc}
\begin{array}{c}
(145, 32, 25, 240)\\
\begin{MYmatrix}
0&1&2&3&2\\
0&0&2&2&2\\
0&0&0&1&0\\
0&1&1&0&0\\
0&1&2&1&0
\end{MYmatrix}
\end{array}
&\begin{array}{c}
(146, 32, 26, 240)\\
\begin{MYmatrix}
0&1&3&1&1\\
0&0&2&1&0\\
0&0&0&1&0\\
0&1&2&0&0\\
1&2&2&1&0
\end{MYmatrix}
\end{array}
&\begin{array}{c}
(147, 32, 26, 120)\\
\begin{MYmatrix}
0&1&3&3&2\\
0&0&2&3&1\\
0&0&0&1&1\\
0&0&0&0&0\\
0&1&1&2&0
\end{MYmatrix}
\end{array}
&\begin{array}{c}
(148, 32, 26, 120)\\
\begin{MYmatrix}
0&1&2&2&1\\
0&0&1&1&1\\
0&0&0&0&0\\
0&1&1&0&1\\
0&1&2&2&0
\end{MYmatrix}
\end{array}
&\begin{array}{c}
(149, 31, 23, 240)\\
\begin{MYmatrix}
0&1&2&1&1\\
0&0&2&1&1\\
0&0&0&1&0\\
0&1&1&0&0\\
1&1&2&1&0
\end{MYmatrix}
\end{array}
&\begin{array}{c}
(150, 31, 23, 240)\\
\begin{MYmatrix}
0&1&3&2&2\\
0&0&3&2&1\\
0&0&0&2&0\\
0&1&1&0&0\\
2&2&2&2&0
\end{MYmatrix}
\end{array}
&\begin{array}{c}
(151, 31, 25, 240)\\
\begin{MYmatrix}
0&1&2&2&1\\
0&0&1&2&1\\
0&0&0&1&0\\
0&0&1&0&0\\
1&1&1&1&0
\end{MYmatrix}
\end{array}
&\begin{array}{c}
(152, 31, 25, 240)\\
\begin{MYmatrix}
0&1&2&2&2\\
0&0&1&2&1\\
0&0&0&1&0\\
0&0&1&0&0\\
1&1&1&1&0
\end{MYmatrix}
\end{array}
&\begin{array}{c}
(153, 31, 22, 240)\\
\begin{MYmatrix}
0&1&2&2&1\\
0&0&1&1&1\\
0&0&0&1&0\\
0&1&1&0&1\\
1&1&1&1&0
\end{MYmatrix}
\end{array}
\end{MYarray}
\end{equation*}
\begin{equation*}
\begin{MYarray}{ccccccccc}
\begin{array}{c}
(154, 31, 22, 120)\\
\begin{MYmatrix}
0&1&1&1&1\\
0&0&1&1&0\\
0&1&0&1&0\\
0&1&1&0&0\\
1&1&2&2&0
\end{MYmatrix}
\end{array}
&\begin{array}{c}
(155, 31, 26, 120)\\
\begin{MYmatrix}
0&2&2&3&3\\
0&0&2&2&1\\
0&0&0&2&1\\
0&0&0&0&0\\
0&1&1&3&0
\end{MYmatrix}
\end{array}
&\begin{array}{c}
(156, 30, 24, 240)\\
\begin{MYmatrix}
0&1&3&2&1\\
0&0&2&1&1\\
0&0&0&1&0\\
0&1&1&0&0\\
1&2&2&1&0
\end{MYmatrix}
\end{array}
&\begin{array}{c}
(157, 30, 24, 240)\\
\begin{MYmatrix}
0&2&2&2&1\\
0&0&2&2&1\\
0&0&0&2&0\\
0&0&1&0&1\\
1&1&3&2&0
\end{MYmatrix}
\end{array}
&\begin{array}{c}
(158, 30, 24, 120)\\
\begin{MYmatrix}
0&1&2&2&2\\
0&0&1&1&1\\
0&0&0&1&0\\
0&1&2&0&0\\
1&1&2&1&0
\end{MYmatrix}
\end{array}
&\begin{array}{c}
(159, 30, 24, 120)\\
\begin{MYmatrix}
0&1&2&2&2\\
0&0&1&2&1\\
0&0&0&1&0\\
0&1&2&0&0\\
1&2&2&1&0
\end{MYmatrix}
\end{array}
&\begin{array}{c}
(160, 29, 24, 240)\\
\begin{MYmatrix}
0&1&2&3&1\\
0&0&1&2&0\\
0&1&0&1&1\\
0&1&1&0&0\\
1&1&1&2&0
\end{MYmatrix}
\end{array}
&\begin{array}{c}
(161, 29, 23, 240)\\
\begin{MYmatrix}
0&1&2&2&2\\
0&0&2&2&1\\
0&0&0&1&0\\
0&1&1&0&0\\
0&1&2&1&0
\end{MYmatrix}
\end{array}
&\begin{array}{c}
(162, 29, 23, 240)\\
\begin{MYmatrix}
0&1&2&1&1\\
0&0&2&1&0\\
0&0&0&1&0\\
0&1&1&0&0\\
1&2&2&2&0
\end{MYmatrix}
\end{array}
\end{MYarray}
\end{equation*}
\begin{equation*}
\begin{MYarray}{ccccccccc}
\begin{array}{c}
(163, 29, 23, 240)\\
\begin{MYmatrix}
0&1&2&2&1\\
0&0&1&1&1\\
0&0&0&1&0\\
0&0&1&0&0\\
1&1&1&2&0
\end{MYmatrix}
\end{array}
&\begin{array}{c}
(164, 29, 21, 240)\\
\begin{MYmatrix}
0&1&2&2&1\\
0&0&2&1&1\\
0&0&0&1&0\\
0&1&1&0&0\\
1&1&2&1&0
\end{MYmatrix}
\end{array}
&\begin{array}{c}
(165, 29, 23, 240)\\
\begin{MYmatrix}
0&1&2&2&1\\
0&0&1&2&0\\
0&1&0&2&1\\
0&1&0&0&0\\
1&1&1&2&0
\end{MYmatrix}
\end{array}
&\begin{array}{c}
(166, 29, 24, 240)\\
\begin{MYmatrix}
0&1&2&1&1\\
0&0&1&1&1\\
0&0&0&1&1\\
0&1&1&0&0\\
1&1&2&1&0
\end{MYmatrix}
\end{array}
&\begin{array}{c}
(167, 29, 24, 240)\\
\begin{MYmatrix}
0&1&2&2&1\\
0&0&1&1&0\\
0&1&0&2&1\\
0&1&1&0&0\\
1&1&1&1&0
\end{MYmatrix}
\end{array}
&\begin{array}{c}
(168, 29, 23, 240)\\
\begin{MYmatrix}
0&1&1&2&1\\
0&0&1&1&1\\
0&1&0&1&0\\
1&1&1&0&0\\
1&1&2&1&0
\end{MYmatrix}
\end{array}
&\begin{array}{c}
(169, 29, 21, 120)\\
\begin{MYmatrix}
0&2&2&2&1\\
0&0&2&1&1\\
0&0&0&1&0\\
0&1&1&0&0\\
1&1&2&1&0
\end{MYmatrix}
\end{array}
&\begin{array}{c}
(170, 29, 24, 120)\\
\begin{MYmatrix}
0&1&2&1&1\\
0&0&1&1&1\\
0&0&0&1&0\\
0&1&1&0&1\\
1&1&1&2&0
\end{MYmatrix}
\end{array}
&\begin{array}{c}
(171, 29, 24, 120)\\
\begin{MYmatrix}
0&1&2&3&2\\
0&0&2&2&1\\
0&0&0&1&1\\
0&0&0&0&0\\
0&1&1&2&0
\end{MYmatrix}
\end{array}
\end{MYarray}
\end{equation*}
\begin{equation*}
\begin{MYarray}{ccccccccc}
\begin{array}{c}
(172, 29, 25, 40)\\
\begin{MYmatrix}
0&2&2&2&2\\
0&0&2&2&1\\
0&0&0&0&0\\
0&1&2&0&2\\
0&2&2&1&0
\end{MYmatrix}
\end{array}
&\begin{array}{c}
(173, 28, 24, 240)\\
\begin{MYmatrix}
0&1&2&2&2\\
0&0&1&1&1\\
0&0&0&1&0\\
0&1&1&0&0\\
1&1&1&1&0
\end{MYmatrix}
\end{array}
&\begin{array}{c}
(174, 28, 24, 240)\\
\begin{MYmatrix}
0&1&2&1&2\\
0&0&1&1&2\\
0&0&0&1&1\\
0&1&1&0&2\\
0&1&1&0&0
\end{MYmatrix}
\end{array}
&\begin{array}{c}
(175, 27, 23, 120)\\
\begin{MYmatrix}
0&1&2&2&2\\
0&0&1&2&2\\
0&1&0&2&1\\
0&1&0&0&0\\
1&1&1&1&0
\end{MYmatrix}
\end{array}
&\begin{array}{c}
(176, 27, 23, 120)\\
\begin{MYmatrix}
0&1&2&2&2\\
0&0&2&2&1\\
0&0&0&0&0\\
0&1&1&0&1\\
0&1&2&1&0
\end{MYmatrix}
\end{array}
&\begin{array}{c}
(177, 26, 23, 240)\\
\begin{MYmatrix}
0&2&2&2&1\\
0&0&1&2&1\\
0&0&0&1&1\\
0&1&2&0&0\\
1&1&2&1&0
\end{MYmatrix}
\end{array}
&\begin{array}{c}
(178, 26, 22, 240)\\
\begin{MYmatrix}
0&1&2&2&1\\
0&0&1&1&0\\
0&1&0&1&0\\
0&1&1&0&1\\
1&1&1&1&0
\end{MYmatrix}
\end{array}
&\begin{array}{c}
(179, 26, 23, 240)\\
\begin{MYmatrix}
0&1&2&3&2\\
0&0&2&2&1\\
0&0&0&1&1\\
0&0&1&0&0\\
0&1&1&1&0
\end{MYmatrix}
\end{array}
&\begin{array}{c}
(180, 26, 24, 240)\\
\begin{MYmatrix}
0&1&2&3&2\\
0&0&1&2&1\\
0&0&0&1&0\\
0&1&1&0&0\\
1&1&2&1&0
\end{MYmatrix}
\end{array}
\end{MYarray}
\end{equation*}
\begin{equation*}
\begin{MYarray}{ccccccccc}
\begin{array}{c}
(181, 26, 24, 240)\\
\begin{MYmatrix}
0&1&2&2&1\\
0&0&1&2&1\\
0&0&0&1&0\\
0&1&1&0&1\\
0&1&2&1&0
\end{MYmatrix}
\end{array}
&\begin{array}{c}
(182, 26, 23, 240)\\
\begin{MYmatrix}
0&1&2&3&2\\
0&0&2&2&1\\
0&0&0&1&1\\
0&1&1&0&0\\
1&1&1&1&0
\end{MYmatrix}
\end{array}
&\begin{array}{c}
(183, 26, 23, 240)\\
\begin{MYmatrix}
0&2&2&2&2\\
0&0&1&2&1\\
0&0&0&1&1\\
0&0&1&0&0\\
0&1&2&1&0
\end{MYmatrix}
\end{array}
&\begin{array}{c}
(184, 26, 24, 120)\\
\begin{MYmatrix}
0&1&1&2&1\\
0&0&1&1&1\\
0&0&0&1&0\\
0&1&1&0&0\\
1&1&2&2&0
\end{MYmatrix}
\end{array}
&\begin{array}{c}
(185, 26, 22, 60)\\
\begin{MYmatrix}
0&2&3&3&2\\
0&0&2&1&2\\
0&0&0&0&0\\
0&1&3&0&1\\
0&2&2&1&0
\end{MYmatrix}
\end{array}
&\begin{array}{c}
(186, 26, 22, 40)\\
\begin{MYmatrix}
0&2&4&2&2\\
0&0&2&2&1\\
0&0&0&0&0\\
0&1&2&0&2\\
0&2&2&1&0
\end{MYmatrix}
\end{array}
&\begin{array}{c}
(187, 25, 23, 240)\\
\begin{MYmatrix}
0&1&2&2&1\\
0&0&1&2&1\\
0&0&0&1&0\\
0&1&1&0&1\\
1&1&1&1&0
\end{MYmatrix}
\end{array}
&\begin{array}{c}
(188, 25, 23, 240)\\
\begin{MYmatrix}
0&1&3&2&2\\
0&0&2&1&2\\
0&0&0&1&0\\
0&1&1&0&1\\
1&1&1&1&0
\end{MYmatrix}
\end{array}
&\begin{array}{c}
(189, 25, 23, 240)\\
\begin{MYmatrix}
0&1&2&1&2\\
0&0&1&1&1\\
0&0&0&1&0\\
0&1&1&0&1\\
1&1&2&1&0
\end{MYmatrix}
\end{array}
\end{MYarray}
\end{equation*}
\begin{equation*}
\begin{MYarray}{ccccccccc}
\begin{array}{c}
(190, 25, 22, 240)\\
\begin{MYmatrix}
0&1&2&2&1\\
0&0&1&2&0\\
0&1&0&1&1\\
0&1&1&0&0\\
1&1&1&2&0
\end{MYmatrix}
\end{array}
&\begin{array}{c}
(191, 25, 23, 240)\\
\begin{MYmatrix}
0&1&1&2&1\\
0&0&1&1&0\\
0&1&0&1&1\\
0&1&1&0&0\\
1&1&2&1&0
\end{MYmatrix}
\end{array}
&\begin{array}{c}
(192, 24, 23, 240)\\
\begin{MYmatrix}
0&1&3&2&1\\
0&0&3&1&1\\
0&0&0&1&0\\
0&1&2&0&0\\
1&2&2&1&0
\end{MYmatrix}
\end{array}
&\begin{array}{c}
(193, 24, 23, 240)\\
\begin{MYmatrix}
0&1&2&2&2\\
0&0&2&1&1\\
0&0&0&1&0\\
0&1&2&0&0\\
1&1&3&1&0
\end{MYmatrix}
\end{array}
&\begin{array}{c}
(194, 24, 23, 240)\\
\begin{MYmatrix}
0&1&2&3&2\\
0&0&2&2&2\\
0&0&0&1&0\\
0&1&2&0&0\\
1&2&2&1&0
\end{MYmatrix}
\end{array}
&\begin{array}{c}
(195, 23, 22, 240)\\
\begin{MYmatrix}
0&1&2&2&1\\
0&0&2&1&1\\
0&0&0&1&0\\
0&1&1&0&0\\
1&2&2&1&0
\end{MYmatrix}
\end{array}
&\begin{array}{c}
(196, 23, 22, 240)\\
\begin{MYmatrix}
0&1&4&3&2\\
0&0&3&2&2\\
0&0&0&2&0\\
0&1&1&0&0\\
2&2&2&2&0
\end{MYmatrix}
\end{array}
&\begin{array}{c}
(197, 23, 22, 240)\\
\begin{MYmatrix}
0&1&3&2&2\\
0&0&2&1&1\\
0&0&0&1&0\\
0&1&1&0&0\\
1&1&2&1&0
\end{MYmatrix}
\end{array}
&\begin{array}{c}
(198, 23, 22, 240)\\
\begin{MYmatrix}
0&1&2&2&1\\
0&0&1&2&1\\
0&0&0&1&0\\
0&1&1&0&0\\
1&1&2&1&0
\end{MYmatrix}
\end{array}
\end{MYarray}
\end{equation*}
\begin{equation*}
\begin{MYarray}{ccccccccc}
\begin{array}{c}
(199, 23, 22, 240)\\
\begin{MYmatrix}
0&2&2&3&2\\
0&0&2&1&1\\
0&0&0&1&0\\
0&1&1&0&0\\
1&1&2&1&0
\end{MYmatrix}
\end{array}
&\begin{array}{c}
(200, 23, 22, 240)\\
\begin{MYmatrix}
0&1&1&1&2\\
0&0&1&0&1\\
0&1&0&1&1\\
1&1&1&0&2\\
1&1&1&0&0
\end{MYmatrix}
\end{array}
&\begin{array}{c}
(201, 23, 22, 240)\\
\begin{MYmatrix}
0&1&2&2&2\\
0&0&1&2&1\\
0&0&0&2&0\\
0&1&1&0&0\\
2&2&3&2&0
\end{MYmatrix}
\end{array}
&\begin{array}{c}
(202, 23, 22, 240)\\
\begin{MYmatrix}
0&1&2&2&2\\
0&0&1&2&1\\
0&1&0&1&0\\
1&1&1&0&0\\
1&1&2&1&0
\end{MYmatrix}
\end{array}
&\begin{array}{c}
(203, 23, 22, 120)\\
\begin{MYmatrix}
0&1&2&2&2\\
0&0&1&2&1\\
0&0&0&1&0\\
0&1&2&0&0\\
1&2&2&3&0
\end{MYmatrix}
\end{array}
&\begin{array}{c}
(204, 22, 21, 240)\\
\begin{MYmatrix}
0&1&2&2&1\\
0&0&1&1&0\\
0&1&0&1&1\\
0&1&1&0&0\\
1&1&1&2&0
\end{MYmatrix}
\end{array}
&\begin{array}{c}
(205, 22, 21, 240)\\
\begin{MYmatrix}
0&2&2&2&2\\
0&0&2&1&0\\
0&1&0&2&0\\
0&2&1&0&0\\
2&2&4&3&0
\end{MYmatrix}
\end{array}
&\begin{array}{c}
(206, 22, 21, 240)\\
\begin{MYmatrix}
0&1&2&2&2\\
0&0&2&1&1\\
0&0&0&1&0\\
0&1&1&0&1\\
1&1&1&1&0
\end{MYmatrix}
\end{array}
&\begin{array}{c}
(207, 22, 22, 240)\\
\begin{MYmatrix}
0&2&3&2&1\\
0&0&1&1&1\\
0&1&0&2&0\\
1&1&1&0&0\\
1&1&2&2&0
\end{MYmatrix}
\end{array}
\end{MYarray}
\end{equation*}
\begin{equation*}
\begin{MYarray}{ccccccccc}
\begin{array}{c}
(208, 21, 21, 240)\\
\begin{MYmatrix}
0&1&2&2&1\\
0&0&2&1&1\\
0&0&0&1&0\\
0&1&2&0&0\\
1&2&2&1&0
\end{MYmatrix}
\end{array}
&\begin{array}{c}
(209, 21, 21, 240)\\
\begin{MYmatrix}
0&1&2&2&2\\
0&0&2&2&1\\
0&0&0&1&0\\
0&1&2&0&0\\
1&2&2&1&0
\end{MYmatrix}
\end{array}
&\begin{array}{c}
(210, 21, 21, 240)\\
\begin{MYmatrix}
0&1&3&2&1\\
0&0&2&1&1\\
0&0&0&1&0\\
0&1&2&0&0\\
1&2&2&1&0
\end{MYmatrix}
\end{array}
&\begin{array}{c}
(211, 21, 21, 240)\\
\begin{MYmatrix}
0&1&2&2&2\\
0&0&2&1&1\\
0&0&0&1&0\\
0&1&2&0&0\\
1&1&2&1&0
\end{MYmatrix}
\end{array}
&\begin{array}{c}
(212, 21, 21, 120)\\
\begin{MYmatrix}
0&2&3&2&1\\
0&0&2&1&1\\
0&0&0&1&0\\
0&1&1&0&0\\
1&1&2&1&0
\end{MYmatrix}
\end{array}
&\begin{array}{c}
(213, 20, 20, 240)\\
\begin{MYmatrix}
0&1&2&2&1\\
0&0&1&1&1\\
0&1&0&1&0\\
1&1&1&0&0\\
1&1&2&1&0
\end{MYmatrix}
\end{array}
&\begin{array}{c}
(214, 20, 20, 240)\\
\begin{MYmatrix}
0&2&2&2&2\\
0&0&2&1&1\\
0&0&0&1&0\\
0&1&1&0&0\\
1&1&2&1&0
\end{MYmatrix}
\end{array}
&\begin{array}{c}
(215, 20, 20, 240)\\
\begin{MYmatrix}
0&1&2&2&1\\
0&0&1&2&1\\
0&1&0&1&1\\
0&1&1&0&0\\
1&1&1&1&0
\end{MYmatrix}
\end{array}
&\begin{array}{c}
(216, 19, 19, 240)\\
\begin{MYmatrix}
0&2&2&2&2\\
0&0&2&1&0\\
0&1&0&2&0\\
0&2&1&0&1\\
2&2&2&2&0
\end{MYmatrix}
\end{array}
\end{MYarray}
\end{equation*}
\begin{equation*}
\begin{MYarray}{ccccccccc}
\begin{array}{c}
(217, 19, 19, 240)\\
\begin{MYmatrix}
0&2&3&2&1\\
0&0&1&2&0\\
0&1&0&2&1\\
0&2&1&0&0\\
2&2&2&2&0
\end{MYmatrix}
\end{array}
&\begin{array}{c}
(218, 19, 19, 240)\\
\begin{MYmatrix}
0&2&3&3&1\\
0&0&1&2&0\\
0&1&0&1&1\\
0&2&1&0&0\\
2&2&2&2&0
\end{MYmatrix}
\end{array}
&\begin{array}{c}
(219, 19, 19, 240)\\
\begin{MYmatrix}
0&2&3&2&1\\
0&0&1&2&0\\
0&1&0&1&1\\
0&2&2&0&0\\
2&2&2&2&0
\end{MYmatrix}
\end{array}
&\begin{array}{c}
(220, 19, 19, 240)\\
\begin{MYmatrix}
0&2&3&3&1\\
0&0&1&2&0\\
0&1&0&1&1\\
0&2&2&0&0\\
2&2&2&2&0
\end{MYmatrix}
\end{array}
&\begin{array}{c}
(221, 19, 19, 240)\\
\begin{MYmatrix}
0&2&3&2&2\\
0&0&2&1&0\\
0&1&0&2&1\\
0&2&1&0&0\\
2&2&2&3&0
\end{MYmatrix}
\end{array}
&\begin{array}{c}
(222, 19, 19, 240)\\
\begin{MYmatrix}
0&2&3&2&1\\
0&0&2&2&0\\
0&1&0&2&1\\
0&2&1&0&0\\
2&2&2&2&0
\end{MYmatrix}
\end{array}
&\begin{array}{c}
(223, 19, 19, 240)\\
\begin{MYmatrix}
0&2&3&3&1\\
0&0&1&2&0\\
0&2&0&1&1\\
0&2&2&0&0\\
2&2&2&2&0
\end{MYmatrix}
\end{array}
&\begin{array}{c}
(224, 19, 19, 240)\\
\begin{MYmatrix}
0&2&2&2&2\\
0&0&2&1&0\\
0&1&0&2&0\\
0&2&1&0&0\\
1&2&2&3&0
\end{MYmatrix}
\end{array}
&\begin{array}{c}
(225, 19, 19, 240)\\
\begin{MYmatrix}
0&2&2&2&2\\
0&0&2&1&0\\
0&1&0&2&1\\
0&2&1&0&0\\
2&2&2&3&0
\end{MYmatrix}
\end{array}
\end{MYarray}
\end{equation*}
\begin{equation*}
\begin{MYarray}{cccc}
\begin{array}{c}
(226, 19, 19, 240)\\
\begin{MYmatrix}
0&2&2&3&2\\
0&0&2&1&0\\
0&1&0&2&0\\
0&2&1&0&1\\
2&2&2&2&0
\end{MYmatrix}
\end{array}
&\begin{array}{c}
(227, 19, 19, 240)\\
\begin{MYmatrix}
0&2&2&3&2\\
0&0&2&1&0\\
0&1&0&2&0\\
0&2&1&0&0\\
2&2&3&2&0
\end{MYmatrix}
\end{array}
&\begin{array}{c}
(228, 19, 19, 240)\\
\begin{MYmatrix}
0&2&2&2&2\\
0&0&2&1&0\\
0&1&0&2&0\\
0&2&1&0&0\\
1&2&3&3&0
\end{MYmatrix}
\end{array}
&\begin{array}{c}
(229, 19, 19, 240)\\
\begin{MYmatrix}
0&2&2&3&1\\
0&0&2&1&0\\
0&1&0&1&1\\
1&1&1&0&0\\
1&3&2&2&0
\end{MYmatrix}
\end{array}
\end{MYarray}
\end{equation*}
}


\section{Appendix $2$: Facets of $OMCUT_5$}
Using a modification by the second author of adjacency decomposition method from~\cite{CR}, all facets of $OMCUT_5$ were found: $194$ orbits consisting all together $35320$ facets.

See below the full list of the representatives of orbits presented as $5\times 5$ matrix. The numbers in parenthesis are, respectively, the orbit number, the adjacency of a member of the orbit, the incidence of a member of the orbit, and orbit-size.

{\tiny
\begin{equation*}
\begin{MYarray}{cccccccc}
\begin{array}{c}
(1, 10695, 307, 20)\\
\begin{MYmatrix}
0&0&0&0&0\\
0&0&0&0&0\\
0&0&0&0&0\\
0&0&0&0&1\\
0&0&0&0&0
\end{MYmatrix}
\end{array}
&\begin{array}{c}
(2, 6451, 307, 60)\\
\begin{MYmatrix}
0&0&0&0&0\\
-1&0&0&0&1\\
0&0&0&0&0\\
0&0&0&0&0\\
1&0&0&0&0
\end{MYmatrix}
\end{array}
&\begin{array}{c}
(3, 1590, 198, 60)\\
\begin{MYmatrix}
0&-1&0&1&0\\
-1&0&0&1&0\\
0&0&0&0&0\\
0&0&0&0&1\\
1&1&0&-1&0
\end{MYmatrix}
\end{array}
&\begin{array}{c}
(4, 940, 111, 120)\\
\begin{MYmatrix}
0&-1&0&1&1\\
-1&0&0&1&1\\
0&0&0&0&0\\
0&0&0&0&0\\
1&1&0&-1&0
\end{MYmatrix}
\end{array}
&\begin{array}{c}
(5, 444, 70, 60)\\
\begin{MYmatrix}
0&-1&0&1&1\\
-1&0&0&1&1\\
-1&-1&0&1&1\\
1&1&0&0&-1\\
1&1&0&-1&0
\end{MYmatrix}
\end{array}
&\begin{array}{c}
(6, 352, 111, 120)\\
\begin{MYmatrix}
0&-1&0&0&1\\
-1&0&0&0&1\\
-1&-1&0&1&1\\
1&1&0&0&-1\\
1&1&0&0&0
\end{MYmatrix}
\end{array}
&\begin{array}{c}
(7, 269, 111, 60)\\
\begin{MYmatrix}
0&0&0&0&0\\
-1&0&0&1&1\\
-1&-1&0&1&1\\
1&1&0&0&-1\\
1&1&0&-1&0
\end{MYmatrix}
\end{array}
&\begin{array}{c}
(8, 166, 56, 30)\\
\begin{MYmatrix}
0&-1&1&1&0\\
-1&0&1&1&0\\
0&0&0&-1&1\\
0&0&-1&0&1\\
1&1&0&0&0
\end{MYmatrix}
\end{array}
\end{MYarray}
\end{equation*}
\begin{equation*}
\begin{MYarray}{cccccccc}
\begin{array}{c}
(9, 162, 125, 120)\\
\begin{MYmatrix}
0&-1&0&1&0\\
-1&0&1&0&0\\
0&0&0&1&0\\
0&1&-1&0&1\\
1&0&1&-1&0
\end{MYmatrix}
\end{array}
&\begin{array}{c}
(10, 160, 112, 240)\\
\begin{MYmatrix}
0&-1&1&0&0\\
-1&0&0&1&0\\
-1&0&0&1&1\\
1&0&0&0&0\\
1&1&0&-1&0
\end{MYmatrix}
\end{array}
&\begin{array}{c}
(11, 154, 39, 40)\\
\begin{MYmatrix}
0&-1&-1&1&2\\
-1&0&-1&1&2\\
-1&-1&0&1&2\\
1&1&1&0&-2\\
1&1&1&-1&0
\end{MYmatrix}
\end{array}
&\begin{array}{c}
(12, 129, 61, 240)\\
\begin{MYmatrix}
0&0&0&1&0\\
-1&0&0&1&1\\
-1&0&0&1&0\\
1&0&0&0&0\\
1&0&1&-1&0
\end{MYmatrix}
\end{array}
&\begin{array}{c}
(13, 110, 63, 240)\\
\begin{MYmatrix}
0&-1&0&1&1\\
-1&0&0&1&0\\
-1&0&0&1&1\\
1&0&0&0&0\\
1&1&0&-1&0
\end{MYmatrix}
\end{array}
&\begin{array}{c}
(14, 107, 39, 20)\\
\begin{MYmatrix}
0&-1&-1&2&2\\
-1&0&-1&2&2\\
-1&-1&0&2&2\\
1&1&1&0&-2\\
1&1&1&-2&0
\end{MYmatrix}
\end{array}
&\begin{array}{c}
(15, 99, 81, 120)\\
\begin{MYmatrix}
0&-2&0&1&1\\
-1&0&1&0&0\\
-1&-2&0&2&2\\
1&2&0&0&-1\\
1&2&0&-1&0
\end{MYmatrix}
\end{array}
&\begin{array}{c}
(16, 98, 88, 60)\\
\begin{MYmatrix}
0&-1&1&0&0\\
-1&0&1&0&0\\
-1&-1&0&2&2\\
1&1&0&0&-1\\
1&1&0&-1&0
\end{MYmatrix}
\end{array}
\end{MYarray}
\end{equation*}
\begin{equation*}
\begin{MYarray}{cccccccc}
\begin{array}{c}
(17, 94, 55, 240)\\
\begin{MYmatrix}
0&-1&0&1&0\\
-1&0&1&1&0\\
0&0&0&0&0\\
0&1&-1&0&1\\
1&0&1&0&0
\end{MYmatrix}
\end{array}
&\begin{array}{c}
(18, 92, 61, 240)\\
\begin{MYmatrix}
0&0&-1&1&1\\
-1&0&1&0&1\\
0&0&0&1&0\\
0&0&0&0&0\\
1&0&1&-1&0
\end{MYmatrix}
\end{array}
&\begin{array}{c}
(19, 92, 59, 240)\\
\begin{MYmatrix}
0&-1&1&1&0\\
-1&0&1&0&1\\
0&0&0&0&0\\
0&1&-1&0&1\\
1&0&0&0&0
\end{MYmatrix}
\end{array}
&\begin{array}{c}
(20, 90, 43, 120)\\
\begin{MYmatrix}
0&-2&0&2&2\\
-1&0&0&1&1\\
-1&0&0&1&1\\
1&1&0&0&-1\\
1&1&0&-1&0
\end{MYmatrix}
\end{array}
&\begin{array}{c}
(21, 89, 85, 240)\\
\begin{MYmatrix}
0&-1&0&1&0\\
-1&0&1&2&-1\\
-1&0&0&0&1\\
1&1&-1&0&1\\
1&1&0&-1&0
\end{MYmatrix}
\end{array}
&\begin{array}{c}
(22, 88, 73, 120)\\
\begin{MYmatrix}
0&1&-2&0&1\\
-1&0&0&1&1\\
-2&1&0&-1&2\\
1&0&1&0&-1\\
2&-1&1&1&0
\end{MYmatrix}
\end{array}
&\begin{array}{c}
(23, 86, 43, 80)\\
\begin{MYmatrix}
0&0&-1&1&1\\
-1&0&0&1&1\\
0&-1&0&1&1\\
0&0&0&0&0\\
1&1&1&-1&0
\end{MYmatrix}
\end{array}
&\begin{array}{c}
(24, 83, 74, 240)\\
\begin{MYmatrix}
0&0&1&0&0\\
-1&0&1&1&1\\
-1&-1&0&1&2\\
1&1&0&0&-2\\
2&2&-1&-2&0
\end{MYmatrix}
\end{array}
\end{MYarray}
\end{equation*}
\begin{equation*}
\begin{MYarray}{cccccccc}
\begin{array}{c}
(25, 81, 41, 240)\\
\begin{MYmatrix}
0&0&0&1&1\\
-1&0&0&1&2\\
-1&-1&0&1&2\\
1&1&0&0&-2\\
2&1&0&-2&0
\end{MYmatrix}
\end{array}
&\begin{array}{c}
(26, 80, 51, 120)\\
\begin{MYmatrix}
0&-1&1&0&1\\
-1&0&1&0&1\\
-1&-1&0&2&1\\
1&1&0&0&-1\\
1&1&0&-1&0
\end{MYmatrix}
\end{array}
&\begin{array}{c}
(27, 79, 56, 120)\\
\begin{MYmatrix}
0&1&-1&0&0\\
-1&0&0&1&1\\
-2&0&0&1&1\\
2&0&1&0&-1\\
2&0&1&-1&0
\end{MYmatrix}
\end{array}
&\begin{array}{c}
(28, 78, 51, 240)\\
\begin{MYmatrix}
0&0&0&1&0\\
-1&0&1&0&0\\
0&0&0&1&1\\
0&0&0&0&0\\
1&1&0&-1&0
\end{MYmatrix}
\end{array}
&\begin{array}{c}
(29, 77, 51, 240)\\
\begin{MYmatrix}
0&0&1&1&0\\
-1&0&0&1&1\\
-1&0&0&1&1\\
1&0&0&0&-1\\
2&0&1&-2&0
\end{MYmatrix}
\end{array}
&\begin{array}{c}
(30, 69, 42, 240)\\
\begin{MYmatrix}
0&0&1&0&0\\
-1&0&0&0&1\\
-1&-1&0&2&1\\
1&1&0&0&-1\\
1&1&0&0&0
\end{MYmatrix}
\end{array}
&\begin{array}{c}
(31, 67, 55, 240)\\
\begin{MYmatrix}
0&-1&0&1&1\\
-1&0&1&0&0\\
0&0&0&1&1\\
0&0&0&0&0\\
1&1&0&-1&0
\end{MYmatrix}
\end{array}
&\begin{array}{c}
(32, 67, 38, 120)\\
\begin{MYmatrix}
0&-1&0&0&1\\
-1&0&0&0&1\\
-1&-1&0&1&2\\
1&1&1&0&-2\\
2&2&1&-1&0
\end{MYmatrix}
\end{array}
\end{MYarray}
\end{equation*}
\begin{equation*}
\begin{MYarray}{cccccccc}
\begin{array}{c}
(33, 66, 54, 10)\\
\begin{MYmatrix}
0&-1&-1&1&1\\
-1&0&-1&1&1\\
-1&-1&0&1&1\\
1&1&1&0&-1\\
1&1&1&-1&0
\end{MYmatrix}
\end{array}
&\begin{array}{c}
(34, 65, 44, 120)\\
\begin{MYmatrix}
0&-1&1&1&1\\
-1&0&1&1&1\\
-1&-1&0&1&2\\
1&1&0&0&-2\\
2&2&-1&-2&0
\end{MYmatrix}
\end{array}
&\begin{array}{c}
(35, 65, 39, 240)\\
\begin{MYmatrix}
0&0&0&0&0\\
-1&0&1&1&0\\
0&-1&0&1&2\\
1&1&1&0&-2\\
1&2&1&-2&0
\end{MYmatrix}
\end{array}
&\begin{array}{c}
(36, 64, 51, 240)\\
\begin{MYmatrix}
0&-2&0&2&1\\
-1&0&1&1&0\\
0&0&0&0&0\\
0&1&-1&0&1\\
1&1&1&-1&0
\end{MYmatrix}
\end{array}
&\begin{array}{c}
(37, 62, 39, 120)\\
\begin{MYmatrix}
0&-1&0&2&2\\
-1&0&0&2&2\\
-1&-1&0&1&2\\
1&1&0&0&-2\\
2&2&0&-3&0
\end{MYmatrix}
\end{array}
&\begin{array}{c}
(38, 60, 51, 120)\\
\begin{MYmatrix}
0&0&-1&1&0\\
-1&0&1&1&0\\
0&0&0&0&0\\
0&0&0&0&1\\
1&0&1&0&0
\end{MYmatrix}
\end{array}
&\begin{array}{c}
(39, 60, 40, 120)\\
\begin{MYmatrix}
0&-1&0&1&0\\
-1&0&0&1&0\\
-1&-1&0&1&2\\
1&1&1&0&-1\\
2&2&1&-2&0
\end{MYmatrix}
\end{array}
&\begin{array}{c}
(40, 59, 73, 120)\\
\begin{MYmatrix}
0&-1&2&1&0\\
-1&0&2&1&0\\
-1&-1&0&1&3\\
1&1&0&0&-2\\
3&3&-2&-3&0
\end{MYmatrix}
\end{array}
\end{MYarray}
\end{equation*}
\begin{equation*}
\begin{MYarray}{cccccccc}
\begin{array}{c}
(41, 59, 67, 120)\\
\begin{MYmatrix}
0&-1&0&0&1\\
-1&0&1&-1&1\\
-1&0&0&1&1\\
1&0&1&0&-1\\
1&1&-1&1&0
\end{MYmatrix}
\end{array}
&\begin{array}{c}
(42, 58, 47, 240)\\
\begin{MYmatrix}
0&-1&0&1&1\\
-1&0&1&0&1\\
-1&-1&0&1&1\\
0&1&0&0&0\\
2&1&0&-1&0
\end{MYmatrix}
\end{array}
&\begin{array}{c}
(43, 57, 67, 120)\\
\begin{MYmatrix}
0&-1&1&0&1\\
-1&0&1&0&1\\
-1&-1&0&1&2\\
1&1&0&0&-2\\
2&2&-1&-1&0
\end{MYmatrix}
\end{array}
&\begin{array}{c}
(44, 57, 53, 120)\\
\begin{MYmatrix}
0&-1&1&1&0\\
-1&0&1&1&0\\
0&0&0&0&0\\
0&0&0&0&1\\
1&1&-1&0&0
\end{MYmatrix}
\end{array}
&\begin{array}{c}
(45, 57, 43, 240)\\
\begin{MYmatrix}
0&-1&-1&2&1\\
-1&0&0&1&1\\
-1&-1&0&2&1\\
0&1&0&0&0\\
2&1&1&-2&0
\end{MYmatrix}
\end{array}
&\begin{array}{c}
(46, 56, 42, 240)\\
\begin{MYmatrix}
0&-1&0&0&1\\
-1&0&1&-1&2\\
-1&0&0&1&0\\
1&0&1&0&0\\
1&1&-1&1&0
\end{MYmatrix}
\end{array}
&\begin{array}{c}
(47, 56, 37, 240)\\
\begin{MYmatrix}
0&0&1&-1&0\\
-1&0&0&1&1\\
0&-1&0&1&1\\
0&1&1&0&-1\\
1&0&0&1&0
\end{MYmatrix}
\end{array}
&\begin{array}{c}
(48, 53, 63, 40)\\
\begin{MYmatrix}
0&-1&-1&1&1\\
-1&0&-1&1&1\\
-1&-1&0&1&1\\
0&0&0&0&1\\
2&2&2&-2&0
\end{MYmatrix}
\end{array}
\end{MYarray}
\end{equation*}
\begin{equation*}
\begin{MYarray}{cccccccc}
\begin{array}{c}
(49, 53, 45, 240)\\
\begin{MYmatrix}
0&0&1&0&0\\
-1&0&0&0&1\\
-1&-1&0&2&1\\
1&1&0&0&-1\\
1&1&1&-1&0
\end{MYmatrix}
\end{array}
&\begin{array}{c}
(50, 52, 40, 240)\\
\begin{MYmatrix}
0&-2&0&2&2\\
-1&0&1&1&2\\
-1&0&0&1&0\\
1&1&0&0&-1\\
1&2&-1&-1&0
\end{MYmatrix}
\end{array}
&\begin{array}{c}
(51, 50, 44, 120)\\
\begin{MYmatrix}
0&0&1&0&0\\
-1&0&1&0&1\\
0&-1&0&1&1\\
0&1&0&0&-1\\
1&2&-1&0&0
\end{MYmatrix}
\end{array}
&\begin{array}{c}
(52, 50, 34, 240)\\
\begin{MYmatrix}
0&-1&0&1&1\\
-1&0&-1&2&1\\
-1&0&0&2&1\\
1&0&1&0&-1\\
1&1&1&-2&0
\end{MYmatrix}
\end{array}
&\begin{array}{c}
(53, 50, 27, 120)\\
\begin{MYmatrix}
0&-2&0&1&1\\
-1&0&1&0&0\\
-1&0&0&1&1\\
2&1&1&0&-1\\
2&1&1&-1&0
\end{MYmatrix}
\end{array}
&\begin{array}{c}
(54, 49, 63, 40)\\
\begin{MYmatrix}
0&0&-1&1&0\\
-1&0&0&1&0\\
0&-1&0&1&0\\
0&0&0&0&2\\
1&1&1&-1&0
\end{MYmatrix}
\end{array}
&\begin{array}{c}
(55, 49, 51, 120)\\
\begin{MYmatrix}
0&-1&0&1&1\\
-1&0&1&0&0\\
-1&-1&0&2&2\\
1&1&0&0&-1\\
1&1&0&-1&0
\end{MYmatrix}
\end{array}
&\begin{array}{c}
(56, 49, 47, 120)\\
\begin{MYmatrix}
0&-1&1&1&1\\
-1&0&1&0&1\\
-1&0&0&1&1\\
1&0&0&0&-1\\
2&1&-1&-1&0
\end{MYmatrix}
\end{array}
\end{MYarray}
\end{equation*}
\begin{equation*}
\begin{MYarray}{cccccccc}
\begin{array}{c}
(57, 48, 47, 240)\\
\begin{MYmatrix}
0&-2&1&2&0\\
-1&0&1&1&0\\
0&1&0&0&0\\
0&1&-1&0&1\\
1&0&1&-1&0
\end{MYmatrix}
\end{array}
&\begin{array}{c}
(58, 48, 43, 120)\\
\begin{MYmatrix}
0&-1&1&0&1\\
-1&0&1&0&1\\
-1&-1&0&1&1\\
1&1&-1&0&0\\
1&1&0&0&0
\end{MYmatrix}
\end{array}
&\begin{array}{c}
(59, 47, 32, 240)\\
\begin{MYmatrix}
0&-2&1&1&1\\
-1&0&1&0&0\\
-1&-1&0&1&2\\
1&1&0&0&-1\\
1&2&0&0&0
\end{MYmatrix}
\end{array}
&\begin{array}{c}
(60, 46, 61, 60)\\
\begin{MYmatrix}
0&-1&1&-1&1\\
-1&0&1&-1&1\\
-1&-1&0&2&1\\
1&1&-1&0&0\\
1&1&0&1&0
\end{MYmatrix}
\end{array}
&\begin{array}{c}
(61, 46, 30, 240)\\
\begin{MYmatrix}
0&2&-2&0&1\\
-1&0&0&1&1\\
-2&1&0&0&1\\
1&0&1&0&0\\
2&-1&1&0&0
\end{MYmatrix}
\end{array}
&\begin{array}{c}
(62, 45, 49, 240)\\
\begin{MYmatrix}
0&-1&1&0&0\\
-1&0&0&0&1\\
-1&1&0&1&1\\
1&1&0&0&-1\\
2&0&1&-1&0
\end{MYmatrix}
\end{array}
&\begin{array}{c}
(63, 44, 38, 240)\\
\begin{MYmatrix}
0&-2&1&1&1\\
-1&0&1&0&0\\
0&0&0&1&0\\
0&1&0&0&1\\
1&1&0&-1&0
\end{MYmatrix}
\end{array}
&\begin{array}{c}
(64, 43, 43, 120)\\
\begin{MYmatrix}
0&-2&0&1&1\\
-1&0&0&1&1\\
-1&-2&0&2&2\\
1&2&0&0&-1\\
1&2&0&-1&0
\end{MYmatrix}
\end{array}
\end{MYarray}
\end{equation*}
\begin{equation*}
\begin{MYarray}{cccccccc}
\begin{array}{c}
(65, 42, 44, 240)\\
\begin{MYmatrix}
0&-2&0&1&2\\
-1&0&1&0&1\\
-2&-2&0&2&2\\
1&2&0&0&-1\\
2&2&0&-1&0
\end{MYmatrix}
\end{array}
&\begin{array}{c}
(66, 42, 38, 240)\\
\begin{MYmatrix}
0&0&1&0&0\\
-1&0&0&0&1\\
0&-1&0&1&1\\
0&1&1&0&-1\\
1&1&-1&1&0
\end{MYmatrix}
\end{array}
&\begin{array}{c}
(67, 42, 36, 240)\\
\begin{MYmatrix}
0&-1&0&1&1\\
-1&0&1&0&1\\
-1&0&0&1&1\\
1&0&1&0&-1\\
1&1&-1&0&0
\end{MYmatrix}
\end{array}
&\begin{array}{c}
(68, 42, 36, 240)\\
\begin{MYmatrix}
0&-1&1&0&0\\
-1&0&0&0&1\\
-1&0&0&1&2\\
1&1&1&0&-2\\
2&1&1&-1&0
\end{MYmatrix}
\end{array}
&\begin{array}{c}
(69, 41, 41, 240)\\
\begin{MYmatrix}
0&-2&1&1&1\\
-1&0&0&1&1\\
0&0&0&0&1\\
0&1&1&0&-1\\
1&1&-1&0&0
\end{MYmatrix}
\end{array}
&\begin{array}{c}
(70, 41, 40, 240)\\
\begin{MYmatrix}
0&-1&0&1&1\\
-1&0&1&2&0\\
-1&0&0&0&1\\
1&1&-1&0&0\\
1&1&0&-1&0
\end{MYmatrix}
\end{array}
&\begin{array}{c}
(71, 41, 40, 240)\\
\begin{MYmatrix}
0&0&1&2&0\\
-1&0&0&2&2\\
-1&1&0&1&0\\
1&0&0&0&-1\\
2&0&2&-3&0
\end{MYmatrix}
\end{array}
&\begin{array}{c}
(72, 41, 39, 240)\\
\begin{MYmatrix}
0&0&0&0&0\\
-1&0&1&1&1\\
-1&-1&0&1&2\\
1&1&1&0&-2\\
2&2&0&-2&0
\end{MYmatrix}
\end{array}
\end{MYarray}
\end{equation*}
\begin{equation*}
\begin{MYarray}{cccccccc}
\begin{array}{c}
(73, 41, 37, 240)\\
\begin{MYmatrix}
0&0&1&-1&1\\
-1&0&0&1&1\\
-1&-1&0&1&1\\
1&0&0&0&0\\
1&1&0&0&0
\end{MYmatrix}
\end{array}
&\begin{array}{c}
(74, 41, 35, 240)\\
\begin{MYmatrix}
0&-2&0&2&2\\
-1&0&-1&2&2\\
-1&-1&0&2&1\\
1&1&0&0&-1\\
1&2&1&-2&0
\end{MYmatrix}
\end{array}
&\begin{array}{c}
(75, 40, 44, 240)\\
\begin{MYmatrix}
0&-2&2&0&1\\
-1&0&1&0&0\\
0&0&0&1&0\\
0&1&0&0&1\\
1&1&-1&0&0
\end{MYmatrix}
\end{array}
&\begin{array}{c}
(76, 40, 41, 120)\\
\begin{MYmatrix}
0&-1&2&1&0\\
-1&0&2&1&0\\
-1&-1&0&1&3\\
1&1&0&0&-2\\
2&2&-1&-2&0
\end{MYmatrix}
\end{array}
&\begin{array}{c}
(77, 40, 29, 240)\\
\begin{MYmatrix}
0&1&2&0&0\\
-1&0&1&0&1\\
0&-1&0&2&1\\
1&1&0&0&-2\\
1&2&-1&0&0
\end{MYmatrix}
\end{array}
&\begin{array}{c}
(78, 39, 38, 240)\\
\begin{MYmatrix}
0&-2&0&2&2\\
-1&0&0&1&0\\
-1&1&0&1&2\\
1&0&1&0&-1\\
1&1&0&-1&0
\end{MYmatrix}
\end{array}
&\begin{array}{c}
(79, 39, 37, 120)\\
\begin{MYmatrix}
0&-1&1&0&0\\
-1&0&1&0&0\\
-1&-1&0&2&3\\
1&1&1&0&-2\\
2&2&1&-2&0
\end{MYmatrix}
\end{array}
&\begin{array}{c}
(80, 39, 36, 240)\\
\begin{MYmatrix}
0&0&1&1&0\\
-1&0&1&1&1\\
-1&-1&0&1&2\\
1&1&0&0&-2\\
2&1&0&-2&0
\end{MYmatrix}
\end{array}
\end{MYarray}
\end{equation*}
\begin{equation*}
\begin{MYarray}{cccccccc}
\begin{array}{c}
(81, 39, 36, 240)\\
\begin{MYmatrix}
0&0&1&0&0\\
-1&0&1&0&0\\
0&-1&0&1&2\\
0&1&0&0&0\\
1&1&0&-1&0
\end{MYmatrix}
\end{array}
&\begin{array}{c}
(82, 39, 32, 240)\\
\begin{MYmatrix}
0&0&1&0&0\\
-1&0&1&1&1\\
-1&-1&0&1&2\\
1&1&0&0&-2\\
2&1&0&-1&0
\end{MYmatrix}
\end{array}
&\begin{array}{c}
(83, 39, 31, 240)\\
\begin{MYmatrix}
0&-1&-1&2&2\\
-1&0&0&2&2\\
-1&-1&0&1&2\\
1&1&1&0&-2\\
1&2&0&-2&0
\end{MYmatrix}
\end{array}
&\begin{array}{c}
(84, 38, 38, 120)\\
\begin{MYmatrix}
0&-2&0&2&2\\
-1&0&1&0&0\\
-1&0&0&2&2\\
1&1&0&0&-1\\
1&1&0&-1&0
\end{MYmatrix}
\end{array}
&\begin{array}{c}
(85, 38, 30, 240)\\
\begin{MYmatrix}
0&0&2&-1&1\\
-1&0&1&1&1\\
0&1&0&1&-2\\
1&-1&2&0&1\\
1&1&-3&1&0
\end{MYmatrix}
\end{array}
&\begin{array}{c}
(86, 37, 45, 40)\\
\begin{MYmatrix}
0&-1&-1&2&1\\
-1&0&-1&2&1\\
-1&-1&0&2&1\\
0&0&0&0&1\\
2&2&2&-3&0
\end{MYmatrix}
\end{array}
&\begin{array}{c}
(87, 35, 39, 240)\\
\begin{MYmatrix}
0&-1&0&1&1\\
-1&0&1&1&0\\
0&0&0&0&1\\
0&0&-1&0&1\\
1&2&0&-1&0
\end{MYmatrix}
\end{array}
&\begin{array}{c}
(88, 35, 39, 120)\\
\begin{MYmatrix}
0&-1&0&1&2\\
-1&0&0&1&2\\
-2&-2&0&2&3\\
2&2&0&0&-3\\
2&2&0&-2&0
\end{MYmatrix}
\end{array}
\end{MYarray}
\end{equation*}
\begin{equation*}
\begin{MYarray}{cccccccc}
\begin{array}{c}
(89, 35, 36, 240)\\
\begin{MYmatrix}
0&-1&0&1&1\\
-1&0&1&1&1\\
-2&-2&0&2&2\\
1&2&0&0&-1\\
2&2&0&-2&0
\end{MYmatrix}
\end{array}
&\begin{array}{c}
(90, 35, 33, 240)\\
\begin{MYmatrix}
0&1&-2&1&1\\
-1&0&0&1&1\\
-2&1&0&0&2\\
1&0&1&0&-1\\
2&-1&1&0&0
\end{MYmatrix}
\end{array}
&\begin{array}{c}
(91, 34, 36, 240)\\
\begin{MYmatrix}
0&-1&1&1&0\\
-1&0&1&1&0\\
-1&0&0&1&2\\
2&2&0&0&-2\\
3&2&0&-3&0
\end{MYmatrix}
\end{array}
&\begin{array}{c}
(92, 34, 33, 240)\\
\begin{MYmatrix}
0&0&1&1&0\\
-1&0&0&1&2\\
-1&-1&0&1&2\\
1&1&0&0&-2\\
2&1&0&-2&0
\end{MYmatrix}
\end{array}
&\begin{array}{c}
(93, 34, 32, 240)\\
\begin{MYmatrix}
0&-1&0&2&2\\
-1&0&1&2&1\\
-1&-1&0&1&2\\
1&1&0&0&-2\\
2&2&0&-3&0
\end{MYmatrix}
\end{array}
&\begin{array}{c}
(94, 33, 34, 240)\\
\begin{MYmatrix}
0&-2&1&2&1\\
-1&0&1&1&1\\
0&1&0&0&-1\\
0&1&-1&0&1\\
1&1&0&-1&0
\end{MYmatrix}
\end{array}
&\begin{array}{c}
(95, 32, 36, 240)\\
\begin{MYmatrix}
0&0&1&1&0\\
-1&0&-2&1&3\\
-1&-2&0&1&2\\
1&2&1&0&-2\\
2&0&0&0&0
\end{MYmatrix}
\end{array}
&\begin{array}{c}
(96, 32, 33, 120)\\
\begin{MYmatrix}
0&0&0&1&0\\
-1&0&-1&1&1\\
-1&-1&0&1&1\\
0&1&1&0&0\\
2&2&2&-2&0
\end{MYmatrix}
\end{array}
\end{MYarray}
\end{equation*}
\begin{equation*}
\begin{MYarray}{cccccccc}
\begin{array}{c}
(97, 32, 31, 240)\\
\begin{MYmatrix}
0&-1&1&2&2\\
-1&0&1&2&1\\
-1&0&0&1&1\\
2&1&-1&0&-2\\
2&1&0&-3&0
\end{MYmatrix}
\end{array}
&\begin{array}{c}
(98, 32, 27, 240)\\
\begin{MYmatrix}
0&0&1&1&1\\
-1&0&1&2&1\\
-1&-1&0&2&2\\
2&1&0&0&-3\\
2&2&-1&-3&0
\end{MYmatrix}
\end{array}
&\begin{array}{c}
(99, 31, 41, 240)\\
\begin{MYmatrix}
0&-1&-1&2&1\\
-1&0&1&0&1\\
-1&-1&0&2&1\\
0&1&0&0&0\\
2&1&1&-2&0
\end{MYmatrix}
\end{array}
&\begin{array}{c}
(100, 31, 41, 240)\\
\begin{MYmatrix}
0&0&1&1&0\\
-1&0&1&1&0\\
-1&0&0&1&2\\
1&0&0&0&-1\\
2&2&-1&-2&0
\end{MYmatrix}
\end{array}
&\begin{array}{c}
(101, 31, 37, 240)\\
\begin{MYmatrix}
0&0&0&0&1\\
-1&0&1&0&1\\
-1&0&0&0&1\\
1&-1&1&0&0\\
1&1&-1&1&0
\end{MYmatrix}
\end{array}
&\begin{array}{c}
(102, 31, 36, 120)\\
\begin{MYmatrix}
0&-1&-1&1&1\\
-1&0&-1&1&1\\
-1&-1&0&1&1\\
0&1&1&0&0\\
2&2&2&-2&0
\end{MYmatrix}
\end{array}
&\begin{array}{c}
(103, 31, 35, 240)\\
\begin{MYmatrix}
0&-1&1&0&1\\
-1&0&0&1&1\\
-1&0&0&1&1\\
1&0&1&0&-1\\
1&1&0&-1&0
\end{MYmatrix}
\end{array}%
\end{MYarray}
\end{equation*}
}

\end{document}